\documentclass[10pt, reqno]{amsart}
\usepackage{amssymb}
\usepackage{amsmath,amscd}
\usepackage[initials,nobysame,alphabetic]{amsrefs}
\usepackage{amsfonts}
\usepackage{mathrsfs}
\usepackage{mathpazo}
\usepackage[arrow,matrix,curve,cmtip,ps]{xy}

\usepackage{amsthm}

\usepackage{float}
\usepackage{hyperref}
\usepackage{graphics}
\usepackage{graphicx}
\usepackage{verbatim}
\usepackage{multirow}
\usepackage{tikz}
\usepackage{enumerate}
\usepackage{eufrak}
\allowdisplaybreaks

\newtheorem{theorem}{Theorem}[section]
\newtheorem{lemma}[theorem]{Lemma}

\newtheorem{corollary}[theorem]{Corollary}

\newtheorem*{theorem*}{Theorem}
\theoremstyle{remark}
\newtheorem{remark}[theorem]{Remark}
\newtheorem{definition}[theorem]{Definition}
\newtheorem{example}[theorem]{Example}

\numberwithin{equation}{section}

\newcommand{\R}{\mathbb{R}}

\newcommand{\E}{\mathbb{E}}
\newcommand{\Prob}{\mathbb{P}}

\newcommand{\mytilde}{\raise.17ex\hbox{$\scriptstyle\mathtt{\sim}$}}

\begin{document}
\title[The Principal-Agent Problem]{The Principal-Agent Problem; A Stochastic Maximum Principle Approach}

\author{Boualem Djehiche and Peter Helgesson}

\address{Department of Mathematics \\ KTH Royal Institute of Technology \\ 100 44, Stockholm \\ Sweden}
\email{boualem@math.kth.se}
\address{Department of Mathematics \\ Chalmers University of Technology \\ 412 96, G\"{o}teborg \\ Sweden}
\email{helgessp@chalmers.se}


\date{\today}

\subjclass[2010]{93E20, 49N70, 49N90}

\keywords{Principal-Agent Problem, Stochastic Maximum Principle, Pontryagin's Maximum Principle}

\begin{abstract}

    We study a general class of Principal-Agent problems in continuous time under hidden action. By formulating the model as a coupled stochastic optimal control problem we are able to find a set of necessary conditions characterizing optimal contracts, using the stochastic maximum principle. An example is carried out to illustrate the proposed approach to the Principal-Agent problem under linear stochastic dynamics with a quadratic performance function.
\\

\end{abstract}

\maketitle

\tableofcontents

\section{Introduction}

Providing incentives between parts in a shared economic environment is fundamental in any modern society. Important examples arise for instance in wage negotiations between employer and employee, terms of insurance between insurance company and client and in contract formation between shareholder and portfolio manager. In the economics literature these problems are known as \textit{Principal-Agent problems}, all having the common feature that a principal hires an agent to work in a certain project, possibly under moral-hazard, and wishes to find an optimal way of providing incentives. In this paper we present a general approach for characterizing such optimal contracts when the action of the agent is hidden information (i.e. moral-hazard) and incentives are delivered continuously in time.\\
\indent The precise structure of the
Principal-Agent problem goes as follows: The principal employs the agent to manage a certain well-defined noisy asset over a fixed period of time. For his/her efforts the agent receives a compensation according to some agreement, set before the period starts. It could for instance involve a lump-sum payment at the end of the period, a continuously paying cash-flow during the period, or both. Depending on what information the principle has at hand to form such an agreement, one distinguishes between two distinct cases; the \textit{Full Information}- and the \textit{Hidden Action-problem}. The full information case differs from the hidden action case in that the principal can observe the actions of the agent in addition to the evolution of the asset. Therefore, under full information the principal is allowed to tailor a contract based on both outcome and effort, not only outcome as for hidden actions. In both cases the contract is constrained by the agent via a so called \textit{participation constraint}, clarifying the minimum requirements of the agent to engage in the project. Under hidden action the contract is further constrained by the \textit{incentive compatibility} condition, meaning that as soon as a contract is assigned the agent will act as to maximize his/her own utility and not necessarily that of the principal.\\
\indent The first paper in which a continuous time version of the Principal-Agent problem appears is \cite{MR882097} by Holmstr\"{o}m and Milgrom. They consider a model in which the agent takes action in continuous time over a finite period and gets rewarded by the principal at the end of this period. In particular, for exponential utility functions they prove the optimal contract to be linear with respect to the output. In 1993 Sch\"{a}tter and Sung \cite{MR1252335} generalized this seminal paper by suggesting a general mathematical framework based on the dynamic programming principle and martingale methods to characterize implementable contracts for a rich class of continuous time models, still though involving lump-sum payments.

In recent years the interest in continuous time versions of the Principal-Agent problem has flourished, mainly because of the available mathematical machinery that offers tractable ways to resolve technical difficulties in discrete time models. The literature has therefore grown substantially, with contributions from many authors, including  Cvitani\'{c}, Wan and Zhang \cite{MR2465709}, Sannikov \cite{MR2433118}, Westerfield \cite{Westerfield} and Williams  \cite{Williams}. A thorough presentation of the field (with an emphasis on the mathematical aspects that arise) can be found in the book \cite{MR2963805} by Cvitani\'{c} and Zhang.\\
\indent In the present literature the paper closest to ours is \cite{Williams}. The aim of that article is to characterize optimal contracts as continuously paying cash-flows in a very general setting, both in full information and under hidden action (with the additional possibility to involve hidden savings of the agent). His idea is to attack the problem by applying the stochastic maximum principle suggested in Bismut \cite{MR0469466}. An explicit example is constructed as a fully dynamical analogue of the exponential utility model by Holmstr\"{o}m and Milgrom, and is solved by using methods involving the dynamical programming principle and the HJB-equation. The book \cite{MR2963805} by Cvitani\'{c} and Zhang has a similar approach, covering a very general class of Principal-Agent problems (involving both lump-sum payments and cash-flows) and using techniques from stochastic optimal control, including martingale methods.\\
\indent Our goal is to study a general stochastic and fully dynamical Principal-Agent problem under hidden actions by applying generalizations of Pontryagin's maximum principle. The idea is to first consider the Agent's problem (assuming a given continuously paying cash-flow as an adapted stochastic process) and characterize the optimal effort. By incentive compatibility we then proceed to the Principal's problem which becomes a state constrained optimal control problem for which we characterize cash-flow optimality. Following that scheme we end up with a coupled pair of stochastic optimization problems. The key required to set up the Principal-Agent problem in this way is (as also pointed out by Yong in \cite{MR2645476}) access to fundamental results concerning existence and uniqueness of solutions to forward-backward stochastic differential equations (from now on FBSDEs). The theory of FBSDEs is far from complete but well established in the stochastic analysis literature (see e.g. \cite{MR1901154}, \cite{MR1262970}, \cite{MR1704232}, \cite{MR1675098}, \cite{MR2223916}, \cite{MR1613669}) and constrained optimal control of such equations has been studied for instance in Ji and Zhou \cite{JiZhou} and Ji and Wei \cite{MR3071092}. The reason to approach the problem by means of the stochastic maximum principle is twofold: Firstly, the method is versatile and does not rely on the dynamic programming principle which opens up the possibility to study models where the HJB-equation does not apply, such as non-Markovian dynamics and performance functionals of mean-field type, including risk measures such as the variance. Secondly, it is possible to extend the maximum principle to include state constrained problems without adding any significant difficulty. It is therefore suitable for explicit calculation.\\
\indent Comparing our approach to what is done in \cite{Williams} we are able to state a full characterization of optimality of the Principal-Agent problem (i.e. characterizing optimal controls of both the agent and the principal satisfying the participation constraint) whereas Williams only solves the Agent's problem. We further believe that our method is simpler to overview and more accessible. Our contribution to the existing literature should be regarded as mathematical rather than economical. We present a general framework for solving a class of Principal-Agent problems, without claiming or investigating possible consequences in economy. The main results of our study are presented in Theorem \ref{TheProblemMainThm} and Theorem \ref{HiddenActionMainThm} in which full characterizations of optimal contracts is stated for two different models.\\
\indent The paper is organized as follows: In Section 2 we present mathematical results from stochastic optimal control that are required for our purposes. Section 3 is devoted to formulating the classes of Principal-Agent problems under study and to characterize optimal contracts by means of the results in Section 2. In Section 4 we make the general approach of the previous sections concrete by a fully solved example in the linear-quadratic (LQ)-setting. As a numerical result we find an interesting and rather counterintuitive behavior of the optimal contract, suggesting a "win-win" relation between the principal and the agent. Finally, in Section 5 we briefly discuss what we have chosen to call \textit{A View Towards Hidden Actions in the Strong Formulation}, meaning that the principal proposes a cash-flow to the agent explicitly given as a Markovian (or feedback type) function  of output. By a simple example it becomes clear that this case is much more difficult to handle and it does not quite fit into the present framework.

\section{Preliminaries}
\label{Preliminaries}

Let $T > 0$ be a fixed time horizon and $(\Omega, \mathcal{F}, \mathbb{F}, \mathbb{P})$ be a filtered probability space satisfying the usual conditions, on which a $d$-dimensional Brownian motion $W = \{ W_t \}_{t \geq 0}$ is defined. We let $\mathbb{F}$ be the natural filtration generated by $W$, augmented by all $\mathbb{P}$-null sets $\mathcal{N}_{\mathbb{P}}$, i.e.  $\mathbb{F} = \mathcal{F}_t \bigvee \mathcal{N}_{\mathbb{P}}$ where $\mathcal{F}_t := \sigma(\{ W_s \}: 0 \leq s \leq t)$.\\
\indent Consider the following control system:
\begin{equation*}
\left\{
    \begin{array}{l}
    dx(t) = b(t,x(t),u(t)) dt + \sigma (t,x(t)) dW_t, \hspace{0.5cm} t \in (0,T]\\
    x(0) = x_0
    \end{array}
\right.
\end{equation*}
with a cost functional of the form
\begin{equation}
    \mathcal{J}(u(\cdot)) := \mathbb{E} \left[ \int_0^T f(t,x(t),u(t)) dt + h(x(T)) \right].
    \label{IntroCostfunc1}
\end{equation}
In the most general setting we assume that $b:[0,T] \times \R^n \times U \rightarrow \R^n$, $\sigma : [0,T] \times \R^n \rightarrow \R^{n \times d}$, $f: [0,T] \times \R^n \times U \rightarrow \R$ and $h: \R \rightarrow \R $ where
\begin{equation*}
\left\{
    \begin{array}{l}
    b(t,x,u) :=
        \left(
        \begin{array}{l}
            b^1(t,x,u)\\
            \hspace{0.7cm}\vdots\\
            b^n(t,x,u)
        \end{array}
        \right)\vspace{0.2cm}\\
    \sigma(t,x) := (\sigma^1(t,x), \ldots , \sigma^d(t,x))\vspace{0.2cm}\\
    \sigma^j(t,x) :=
        \left(
        \begin{array}{l}
        \sigma^{1j}(t,x)\\
        \hspace{0.6cm}\vdots \\
        \sigma^{nj}(t,x)
        \end{array}
        \right), 1 \leq j \leq d
    \end{array}
\right.
\end{equation*}
We require the following assumptions to hold:
\begin{description}
    \item[(S1)]$(U,d)$ is a separable metric space.\vspace{0.2cm}
    \item[(S2)]The maps $b$, $\sigma$, $f$ and $h$ are measurable and there exists a constant $L > 0$ and a modulus of continuity $\bar{\omega} : [0, \infty) \rightarrow [0, \infty)$ such that for $\varphi(t,x,u) = b(t,x,u), \sigma(t,x,u), f(t,x,u), h(x)$ it holds that
        \begin{equation*}
        \left\{
            \begin{array}{l}
                |\varphi(t,x,u) - \varphi(t, \hat{x}, \hat{u})| \leq L|x - \hat{x}| + \bar{\omega}(d(u,\hat{u})),\\
                \hspace{2.6cm} \forall t \in [0,T],\hspace{0.2cm} x, \hat{x} \in \R^n, u, \hat{u} \in U\\
                |\varphi(t,0,u)| \leq L, \hspace{0.2cm} \forall (t,u) \in [0,T] \times U \vspace{0.2cm}
            \end{array}
        \right.
        \end{equation*}
    \item[(S3)]The maps $b$, $\sigma$, $f$ and $h$ are $C^1$ in $x$. Moreover there exists a constant $L > 0$ and a modulus of continuity $\bar{\omega} : [0, \infty) \rightarrow [0, \infty)$ such that for $\varphi(t,x,u) = b(t,x,u), \sigma(t,x,u), f(t,x,u), h(x)$ it holds that
        \begin{equation*}
        \left\{
            \begin{array}{l}
                |\varphi_x(t,x,u) - \varphi_x(t, \hat{x}, \hat{u})| \leq L|x - \hat{x}| + \bar{\omega}(d(u,\hat{u})),\\
                \hspace{2.5cm} \forall t \in [0,T],\hspace{0.2cm} x, \hat{x} \in \R^n,\hspace{0.2cm} u, \hat{u} \in U
            \end{array}
        \right.
        \end{equation*}
\end{description}
Finally we define the space of \textit{admissible controls} as
\begin{equation*}
    \mathcal{U}[0,T] := \left\{ u: [0,T] \times \Omega \rightarrow U; \,\,\, u \text{ is } \{ \mathcal{F}_t \}_{t \geq 0} \text{-adapted} \right\}
\end{equation*}
and formulate our optimal control problem:
\begin{description}
    \item[Problem (S)] Minimize (\ref{IntroCostfunc1}) over $\mathcal{U}[0,T]$.
\end{description}
Any $\bar{u}(\cdot) \in \mathcal{U}[0,T]$ satisfying
\begin{equation*}
    \mathcal{J}(\bar{u}(\cdot)) = \smash{\displaystyle\inf_{u(\cdot) \in \mathcal{U}[0,T]}} \mathcal{J}(u(\cdot)) \vspace{0.2cm}
\end{equation*}
is called an \textit{optimal control} and the corresponding $\bar{x}(\cdot)$ is called the \textit{optimal state process}. We will refer to $(\bar{x}(\cdot), \bar{u}(\cdot))$ as an \textit{optimal pair}. We are now ready to state the celebrated stochastic maximum principle, providing necessary conditions for optimality in problem (S).
\begin{theorem}[The Stochastic Maximum Principle]
    Let the regularity conditions in (S1)-(S3) hold and consider an optimal pair $(\bar{x}(\cdot), \bar{u}(\cdot))$ of problem (S). Then there exists a pair of processes $(p(\cdot), q(\cdot)) \in L_{\mathcal{F}}^2(0,T;\R^n) \times (L_{\mathcal{F}}^2(0,T;\R^n))^d$, where $q(\cdot) := (q^1(\cdot),\ldots, q^d(\cdot))$, satisfying the adjoint equation
    \begin{equation}
    \left\{
    \begin{array}{lll}
        dp(t) = -\left\{ b_x(t, \bar{x}(t), \bar{u}(t)) p(t) + \sum_{j=1}^d \sigma_x^j(t,\bar{x}(t))^\top q(t)^j - f_x(t, \bar{x}(t), \bar{u}(t)) \right\} dt \\ \quad\quad\quad\hspace{0.25cm}+ q(t) dW_t, \\
        p(T) = -h_x(\bar{x}(T)),
    \end{array}
    \right.
    \end{equation}
    such that
    \begin{equation}
        \bar{u}(t)= \smash{\displaystyle\arg\max_{u \in U}}\hspace{0.1cm}\mathcal{H}(t, \bar{x}(t), u, p(t), q(t)), \hspace{0.3cm} \text{a.e. } t \in [0,T], \hspace{0.3cm} \Prob\text{-a.s.}
    \end{equation}
    where the Hamiltonian function $\mathcal{H}$ is given by
    \begin{equation}
        \mathcal{H} (t, x, u, p, q) := \langle p , b(t,x,u) \rangle + \text{tr} \left[ q^\top \sigma(t,x,u) \right] - f(t,x,u)
    \end{equation}
    for $(t,x,u,p,q) \in [0,T] \times \R^n \times U \times \R^n \times \R^{n \times d}$.
    \label{IntroThmSMP1}
\end{theorem}
\begin{remark}
It is important to remember that Theorem \ref{IntroThmSMP1} merely states a set of necessary conditions for optimality in (S). It does not claim the existence of such. Existence of stochastic optimal controls (both in the strong and the weak sense) has been a subject of study since the sixties (see e.g. \cite{MR0192946}) and, at least for strong solutions, the results seem to depend a lot upon the statement of the problem. In the weak sense, an account of existence results is  to be found in \cite{MR1696772} (Theorem 5.3, p. 71).
\label{IntroRmkExistence}
\end{remark}
\begin{remark}
    Restricting the space $U$ to be convex allows for a controlled diffusion coefficient $\sigma(t,x,u)$, without changing the conclusion of Theorem \ref{IntroThmSMP1}. In the case of a non-convex control space the stochastic maximum principle with controlled diffusion was proven in \cite{MR1051633} and requires a solution $(P,Q)$ of an additional adjoint BSDE. We choose to leave this most general form as reference in order to keep the presentation clear.
    \label{IntroRemarkThmSMP1}
\end{remark}
As pointed out in Remark \ref{IntroRmkExistence} it is a non-trivial task to prove the existence of an optimal pair $(\bar{x}(\cdot), \bar{u}(\cdot))$ in a general stochastic control model. Under the additional assumptions
\begin{description}
    \item[(S4)] The control domain $U$ is a convex body in $\R^k$. The maps $b$, $\sigma$ and $f$ are locally Lipschitz in $u$ and their derivatives in $x$ are continuous in $(x,u)$,
\end{description}
the following theorem provides sufficient conditions for optimality in (S).
\begin{theorem}[Sufficient Conditions for Optimality]
Assume conditions (S1)-(S4) and let $(\bar{x}(\cdot), \bar{u}(\cdot), p(\cdot), q(\cdot))$ be an admissible 4-tuple. Suppose that $h(\cdot)$ is convex, $\mathcal{H} (t, \cdot, \cdot, p(t), q(t))$ is concave for all $t \in [0,T]$ $\mathbb{P}$-a.s. and
\begin{equation*}
    \bar{u}(t)= \smash{\displaystyle\arg\max_{u \in U}}\hspace{0.1cm}\mathcal{H}(t, \bar{x}(t), u, p(t), q(t)), \hspace{0.3cm} \text{a.e. } t \in [0,T], \hspace{0.3cm} \Prob\text{-a.s.}
\end{equation*}
Then $(\bar{x}(\cdot), \bar{u}(\cdot))$ is an optimal pair for problem (S).
\end{theorem}
\noindent For the sake of simplicity, from now on we will only consider the case $d=1$ and assume the functions $b$, $\sigma$, $f$ to be $\R$-valued. We will also restrict our selves to consider the case in which $U \subseteq \R$, equipped with the usual Euclidean metric.\\
\indent The stochastic maximum principle has since the early days of the subject (in pioneering papers by e.g. Bismut and Bensoussan in \cite{MR0469466} and \cite{MR705931}) developed a lot and do by now apply to a wide range of problems more general than (S) (see for instance \cite{MR1051633}, \cite{MR2784835} \cite{MR2822408}, \cite{risksensitive}, to mention a few). For our purposes it will be necessary to have a refined version of Theorem \ref{IntroThmSMP1}, characterizing optimal controls in a FBSDE-dynamical setting under state constraints. More precisely we wish to consider a stochastically controlled system of the form
\begin{equation}
    \left\{
    \begin{array}{l}
        dx(t) = b(t, x(t), y(t), z(t), u(t))dt + \sigma(t, x(t), y(t), z(t))dW_t\\
        dy(t) = -c(t, x(t), y(t), z(t), u(t))dt + z(t) dW_t\\
        x(0) = x_0, \hspace{0.2cm} y(T) = \varphi(x(T)),
    \end{array}
    \right.
    \label{IntroFBSDEsyst}
\end{equation}
with respect to a cost-functional
\begin{equation}
\mathcal{J}(u(\cdot)) := \mathbb{E} \left[ \int_0^T f(t,x(t),y(t),z(t),u(t)) dt + h(x(T)) + g(y(0)) \right],
    \label{IntroCostfunc2}
\end{equation}
but also with respect to a set of state constraints
\begin{equation}
    \begin{array}{l}
    \E \left[ \int_0^T \mathbf{F}(t, x(t), y(t), z(t), u(t)) dt + \mathbf{H} (x(T)) + \mathbf{G}(y(0)) \right] :=\vspace{0.3cm}\\
    \left(
    \begin{array}{l}
        \E \left[ \int_0^T f^1(t, x(t), y(t), z(t), u(t)) dt + h^1(x(T)) + g^1(y(0)) \right]\\
        \hspace{4cm}\vdots\\
        \E \left[ \int_0^T f^l(t, x(t), y(t), z(t), u(t)) dt + h^l(x(T)) + g^l(y(0)) \right]
    \end{array}
    \right) \in \Lambda,
    \end{array}
    \label{IntroStateConstraints}
\end{equation}
for some closed and convex set $\Lambda \subseteq \R^l$. The optimal control problem is:
\begin{description}
    \item[Problem (SC)] Minimize (\ref{IntroCostfunc2}) subject to the state constraints (\ref{IntroStateConstraints}) over the set $\mathcal{U}[0,T]$.
\end{description}
To get a good maximum principle for (SC) we require some further regularity conditions ensuring solvability of (\ref{IntroFBSDEsyst}). These conditions, denoted by (SC1)-(SC3) are stated below and can be found in \cite{MR2230926}. All coefficient functions, and hence all $x(\cdot), y(\cdot), z(\cdot)$, are considered to be real valued.\\
\indent Let $G \in \R \backslash \{ 0 \}$ and introduce the following notation:
\begin{equation*}
    \mathbf{v} =
    \left(
    \begin{array}{l}
        x\\
        y\\
        z
    \end{array}
    \right), \hspace{0.5cm}
    A(t,\mathbf{v}) =
    \left(
    \begin{array}{l}
        -Gf\\
        \hspace{0.25cm}Gb\\
        \hspace{0.25cm}G\sigma
    \end{array}
    \right) (t, \mathbf{v}).
\end{equation*}
We assume the existence of a nonzero constant $G$ so that:
\begin{equation*}
\text{(SC1):}\hspace{0.35cm}
\left\{
    \begin{array}{l}
        \text{(i) } A(t, \mathbf{v})\text{ is uniformly Lipschitz with respect to } \mathbf{v}\vspace{0.2cm}\\
        \text{(ii) } A(t, \mathbf{v}) \in L_{\mathcal{F}}^2(\Omega; C([0,T]; \mathbb{R})) \times L_{\mathcal{F}}^2(\Omega; C([0,T]; \mathbb{R})) \times L_{\mathcal{F}}^2(0,T;\mathbb{R})\\
        \hspace{0.4cm}\text{ for each } \mathbf{v} \in \R^3\vspace{0.2cm}\\
        \text{(iii) } \varphi(x) \text{ is uniformly Lipschitz with respect to } x\in \R\vspace{0.2cm}\\
        \text{(iv) } \varphi(x) \in L^2(\Omega, \mathcal{F}_T, \mathbb{P};\R) \text{ for } x\in \R \vspace{0.2cm}
    \end{array}
\right.
\end{equation*}
\begin{equation*}\text{(SC2):}
\left\{
    \begin{array}{l}
    \langle A(t, \mathbf{v}) - A(t, \tilde{\mathbf{v}}), \mathbf{v} - \tilde{\mathbf{v}} \rangle \leq - \beta_1 |G(x - \tilde{x})|^2\vspace{0.1cm}\\
    \hspace{5cm}- \beta_2 (|G(y - \tilde{y})|^2 + |G(z - \tilde{z})|^2),\vspace{0.2cm}\\
    \langle \varphi(x) - \varphi(\tilde{x}), L(x - \tilde{x}) \rangle \geq \mu |G(x - \tilde{x})|^2,\vspace{0.2cm}\\
    \text{where } \beta_1, \beta_2, \mu \text{ are nonnegative constants with } \beta_1 + \beta_2 > 0, \beta_2 + \mu > 0. \vspace{0.2cm}
    \end{array}
\right.
\end{equation*}
Conditions (SC1) and (SC2) were first introduced in \cite{MR1675098} to study existence and uniqueness of fully coupled FBSDEs. The monotonicity property (SC2) is sometimes called $G$-monotonicity of $A(t, \mathbf{v})$. We further assume that:
\begin{equation*}
\text{(SC3):}\hspace{0.35cm}
\left\{
    \begin{array}{l}
        \text{(i) } b, \sigma , c, \varphi , f, h, g \text{ are continuously differentiable.} \vspace{0.2cm}\\
        \text{(ii) } \text{The derivatives of } b, \sigma, c, \varphi \text{ are bounded.}\vspace{0.2cm}\\
        \text{(iii) } \text{The derivatives of } f \text{ are bounded by } C(1 + |x| + |y| + |z| + |u|).\vspace{0.2cm}\\
        \text{(iv) } \text{The derivatives of } h \text{ and } g \text{ are bounded by } C(1 + |x|) \vspace{0.1cm} \\
        \hspace{0.6cm}\text{ and } C(1 + |y|),\text{ respectively.}\vspace{0.2cm}
    \end{array}
\right.
\end{equation*}

\medskip For convenience, we will use the following notation throughout the paper. For $\phi\in \{b, \sigma,c, f, h\}$, etc, we define,
\begin{equation}\label{notation}\left\{\begin{array}{llll}
\delta\phi(t)=\phi(t,\bar{x}(t), \bar{y}(t), \bar{z}(t),u(t))-\phi(t,\bar{x}(t), \bar{y}(t), \bar{z}(t), \bar{u}(t));\\
\phi_w(t)=\frac{\partial\phi}{\partial w}(t,\bar{x}(t), \bar{y}(t), \bar{z}(t), \bar{u}(t)),\,\,\,w=x,y,z; \\ \phi_{xx}(t)=\frac{\partial^ 2\phi}{\partial x^ 2}(t,\bar{x}(t), \bar{y}(t), \bar{z}(t), \bar{u}(t)).
\end{array}\right.
\end{equation}
whenever, $(\bar{x}(\cdot), \bar{y}(\cdot), \bar{z}(\cdot), \bar{u}(\cdot))$ is an optimal 4-tuple of problem (SC), where $u$ is an admissible control. \\

We are now ready to formulate the state constrained stochastic maximum principle for fully coupled FBSDEs.
\begin{theorem}[The State Constrained Maximum Principle for a FBSDE-system]
    Let the regularity conditions in (SC1)-(SC3) hold and assume $\Lambda \subseteq \R^l$ to be a closed and convex set. If $(\bar{x}(\cdot), \bar{y}(\cdot), \bar{z}(\cdot), \bar{u}(\cdot))$ is an optimal 4-tuple of problem (SC), then there exists a vector $(\lambda_0, \lambda) \in \R^{1+l}$ such that
    \begin{equation}
        \lambda_0 \geq 0, \hspace{0.5cm} |\lambda_0|^2 + |\lambda|^2 = 1,
        \label{IntroLagMult}
    \end{equation}
    satisfying the transversality condition
    \begin{equation}
        \langle \lambda , v - \E \left[ \int_0^T \mathbf{F}(t, {\bar x}(t), {\bar y}(t), {\bar z}(t), {\bar u}(t)) dt + \mathbf{H}({\bar x}(T)) + \mathbf{G}({\bar y}(0)) \right] \rangle \geq 0, \hspace{0.5cm} \forall v \in \Lambda
        \label{IntroTransversalityCond}
    \end{equation}
    and a 3-tuple $(r(\cdot), p(\cdot), q(\cdot)) \in L_{\mathcal{F}}^2(\Omega; C([0,T]; \mathbb{R})) \times L_{\mathcal{F}}^2(\Omega; C([0,T]; \mathbb{R})) \times L_{\mathcal{F}}^2(0,T;\mathbb{R})$ of solutions to the adjoint FBSDE
    \begin{equation}
    \left\{
        \begin{array}{l}
        dr(t)= \left\{ c_y(t) r(t) - b_y(t) p(t) - \sigma_y(t) q(t) + \sum_{i=0}^l \lambda_i f_y^i(t) \right\} dt +\\
        \vspace{0.2cm}\hspace{2.35cm}+\left\{ c_z(t) r(t) - b_z (t)p(t) - \sigma_z(t) q(t) + \sum_{i=0}^l \lambda_i f_z^i(t) \right\} dW_t,\\
        \vspace{0.2cm}dp(t) = -\left\{ -c_x(t) r(t)+ b_x(t) p(t) + \sigma_x (t)q(t) - \sum_{i=0}^l \lambda_i f_x^i (t)\right\} dt + q(t) dW_t,\\
        r(0) = \sum_{i=0}^l \lambda_i \E [ g^i({\bar y}(0)) ], p(T) = -\varphi_x({\bar x}(T))r(T) - \sum_{i=0}^l \lambda_i h_x^i({\bar x}(T)),
        \end{array}
        \label{IntroFBSDEConstSMP}
    \right.
    \end{equation}
    such that
    \begin{equation*}
        \bar{u}(t)= \smash{\displaystyle\arg\max_{u \in U}}\hspace{0.1cm} \mathcal{H}(t, \bar{x}(t),\bar{y}(t),\bar{z}(t),u,r(t),p(t),q(t),\lambda_0,\lambda) \hspace{0.3cm} \text{a.e. } t \in [0,T], \hspace{0.3cm} \Prob\text{-a.s.}
    \end{equation*}
    where the Hamiltonian function $\mathcal{H}$ is given by
    \begin{equation*}
        \begin{array}{l}
        \mathcal{H}(t,x,y,z,u,r,p,q,\lambda_0,\lambda) :=\\
         \hspace{0.5cm}\langle r , -c(t,x,y,z,u) \rangle + \langle p , b(t,x,y,z,u) \rangle + q\sigma(t,x,y,z) - \sum_{i=0}^l \lambda_i f^i(t,x,y,z,u).
        \end{array}
    \end{equation*}
    \label{IntroConstrainedSMPthm}
\end{theorem}
A maximum principle of FBSDEs similar to the above has been studied in \cite{JiZhou} and \cite{MR3071092}, though under so called terminal state constraints. To the best of our knowledge the proof of Theorem \ref{IntroConstrainedSMPthm} is not in the literature, but the arguments closely follow those in the proof of an analogue theorem for controlling SDEs rather than FBSDEs. For the sake of completeness we present the full proof, yet pointing out the close overlap with a similar theorem in \cite{MR1696772} (Theorem 6.1, p. 144). Some necessary but technical results that are required are presented in Appendix \ref{AppendixTechThms}.
\begin{proof}
First we introduce the notation:
\begin{equation*}
    \mathbf{x}(T) :=  \int_0^T \mathbf{F}(t, x(t), y(t), z(t), u(t)) dt + \mathbf{H} (x(T)) + \mathbf{G}(y(0)).
\end{equation*}
By adding a proper constant to the cost functional we may without loss of generality assume that $\mathcal{J}(\bar{u}(\cdot)) = 0$. For any $\rho > 0$ define
\begin{equation*}
    \mathcal{J}_{\rho} (\bar{u}(\cdot)) := \{ [ ( \mathcal{J}(\bar{u}(\cdot)) + \rho )^+ ]^2 + d_{\Lambda} (\mathbb{E}\mathbf{x}(T))^2 \}^{1/2}
\end{equation*}
where $d_{\Lambda}$ is the metric measuring the distance from any point in $\mathbb{R}^l$ to the closed and convex set $\Lambda$, i.e.
\begin{equation*}
    d_{\Lambda}(v) := \smash{\displaystyle\inf_{v' \in \Lambda}}\hspace{0.1cm} |v-v'|, \hspace{0.3cm} \forall v \in \mathbb{R}^l.
\end{equation*}
By the fact that $(\mathcal{U}[0,T],\bar{d})$ is a complete metric space when
\begin{equation*}
    \bar{d}(u(\cdot), \tilde{u}(\cdot)) := |\{ (t,\omega) \in [0,T] \times \Omega: u(t,\omega) \neq \tilde{u}(t,\omega) \}|, \hspace{0.3cm} \forall u(\cdot), \tilde{u}(\cdot) \in \mathcal{U}[0,T],
\end{equation*}
(see Lemma \ref{AppendixMetricLemma}) Ekeland's variational principle (Corollary \ref{AppendixEkelandCor}) applied to $\mathcal{J}_{\rho}$ ensures the existence of a $u_{\rho} (\cdot) \in \mathcal{U}[0,T]$ being optimal to the problem of controlling (\ref{IntroFBSDEsyst}) with respect to the cost functional $\mathcal{J}_{\rho}(u(\cdot)) + \sqrt{\rho} \bar{d}(u_{\rho} (\cdot), u(\cdot))$. Furthermore it holds (also by Ekeland's principle) that $\bar{d}(u_{\rho} (\cdot), \bar{u}(\cdot)) \leq \rho$ so the control $u_{\rho} (\cdot)$ is close to the state constrained optimal control $\bar{u}(\cdot)$. The idea is now to derive necessary conditions for optimality of $u_{\rho} (\cdot)$ and then let $\rho \rightarrow 0$.\\
\indent Fix a $\rho > 0$ and a $u(\cdot) \in \mathcal{U}[0,T]$ and consider the spike variation
\begin{equation*}
    u_{\rho}^{\varepsilon} (t) :=
    \left\{
    \begin{array}{l}
        u_{\rho}(t),\,\, t \in [0,T] \backslash E_{\varepsilon},\\
        u(t),\,\, t \in E_{\varepsilon},
    \end{array}
    \right.
\end{equation*}
for any given $\varepsilon > 0$ and Borel measurable set $E_{\varepsilon} \subseteq [0,T]$ with $|E_{\varepsilon}| = \varepsilon$. Clearly $\bar{d}(u_{\rho}^{\varepsilon} (\cdot), u_{\rho}(\cdot)) \leq |E_{\varepsilon} \times \Omega| = \varepsilon$ so by the optimality of $u_{\rho}(\cdot)$ and a Taylor expansion we get
\begin{equation}
    \begin{array}{l}
    -\sqrt{\rho} \varepsilon \leq \mathcal{J}_{\rho}(u_{\rho}^{\varepsilon}(\cdot)) - \mathcal{J}_{\rho}(u_{\rho}(\cdot)) = \\
    \hspace{0.5cm}\vspace{0.2cm}= \lambda_{\rho}^{0,\varepsilon} (\mathcal{J}(u_{\rho}^{\varepsilon}(\cdot)) - \mathcal{J}(u_{\rho}(\cdot))) + \langle \lambda_{\rho}^{\varepsilon}, \mathbb{E}[\mathbf{x}_{\rho}^{\varepsilon}(T) - \mathbf{x}_{\rho}(T)] \rangle =\\
    \hspace{0.5cm}\vspace{0.2cm}= \mathbb{E}\left[ \sum_{i = 0}^l \lambda_{\rho}^{i,\varepsilon} \left\{ \int_0^T f^i(t, x_{\rho}^{\varepsilon}(t), y_{\rho}^{\varepsilon}(t), z_{\rho}^{\varepsilon}(t), u_{\rho}^{\varepsilon}(t)) - f^i(t, x_{\rho}(t), y_{\rho}(t), z_{\rho}(t), u_{\rho}(t)) dt \right. \right.\\
    \hspace{1.5cm}\vspace{0.2cm} \left. \left. + (h^i(x_{\rho}^{\varepsilon}(T)) - h^i(x_{\rho}(T))) + (g^i(y_{\rho}^{\varepsilon}(0)) - g^i(y_{\rho}(0))) \right\} \right]
    \end{array}
    \label{IntroCSMPestimate1}
\end{equation}
where the linear approximation factors are
\begin{equation}
    \lambda_{\rho}^{0,\varepsilon} := \frac{[\mathcal{J}(u_{\rho}(\cdot)) + \rho]^+}{\mathcal{J}_{\rho}(u_{\rho}(\cdot))} + o(1),
    \label{IntroFactor1}
\end{equation}
\begin{equation}
    \lambda_{\rho}^{\varepsilon} := \frac{d_{\Lambda}( \mathbb{E}[\mathbf{x}_{\rho}(T)] ) \partial d_{\Lambda}( \mathbb{E}[\mathbf{x}_{\rho}(T)] )}{ \mathcal{J}_{\rho}(u_{\rho}(\cdot)) } + o(1),
    \label{IntroFactor2}
\end{equation}
and $\lambda_{\rho}^{\varepsilon} = (\lambda_{\rho}^{1,\varepsilon}, ..., \lambda_{\rho}^{l,\varepsilon}  )$. By the variational inequality in \cite{MR2230926} (Lemma 3.5, p. 163) we get that the inequality in (\ref{IntroCSMPestimate1}) can be estimated further so that
\begin{equation*}
    \begin{array}{l}
    -\sqrt{\rho} \varepsilon \leq \mathbb{E} \left[ \sum_{i = 0}^l \lambda_{\rho}^{i, \varepsilon} \left\{ \int_0^T f_x^i(t)\tilde{x}_{\rho}^{\varepsilon}(t) + f_y^i(t)\tilde{y}_{\rho}^{\varepsilon}(t) + f_z^i(t)\tilde{z}_{\rho}^{\varepsilon}(t) + \delta f^i(t) \chi_{E_{\varepsilon}}(t) dt \right. \right.\\
    \hspace{2.5cm}\vspace{0.2cm} \left. \left. + h_x^i(x_{\rho}(T))\tilde{x}_{\rho}^{\varepsilon}(T) + g_y^i(y_{\rho}(0))\tilde{y}_{\rho}^{\varepsilon}(0) \right\} \right] + o(\varepsilon),
    \end{array}
\end{equation*}
where the processes $\tilde{x}_{\rho}^{\varepsilon}, \tilde{y}_{\rho}^{\varepsilon}, \tilde{z}_{\rho}^{\varepsilon}$ solve the variational FBSDE-system:
\begin{equation*}
\left\{
    \begin{array}{l}
    d\tilde{x}_{\rho}^{\varepsilon}(t) = (b_x(t)\tilde{x}_{\rho}^{\varepsilon}(t) + b_y(t)\tilde{y}_{\rho}^{\varepsilon}(t) + b_z(t)\tilde{z}_{\rho}^{\varepsilon}(t) + \delta b(t)\chi_{E_{\varepsilon}}(t) )dt\\
    \hspace{5cm}\vspace{0.2cm}+ (\sigma_x(t)\tilde{x}_{\rho}^{\varepsilon}(t) + \sigma_y(t)\tilde{y}_{\rho}^{\varepsilon}(t) + \sigma_z(t)\tilde{z}_{\rho}^{\varepsilon}(t)) dW_t,\\
    \vspace{0.2cm} d\tilde{y}_{\rho}^{\varepsilon}(t) = -(f_x(t)\tilde{x}_{\rho}^{\varepsilon}(t) + f_y(t)\tilde{y}_{\rho}^{\varepsilon}(t) + f_z(t)\tilde{z}_{\rho}^{\varepsilon}(t) + \delta f(t)\chi_{E_{\varepsilon}}(t) )dt + \tilde{z}_{\rho}^{\varepsilon}(t) dW_t\\
    \vspace{0.2cm} \tilde{x}_{\rho}^{\varepsilon}(0) = 0, \hspace{0.5cm} \tilde{y}_{\rho}^{\varepsilon}(T) = \varphi_x(x_{\rho}^{\varepsilon}(T))\tilde{x}_{\rho}^{\varepsilon}(T).
    \end{array}
\right.
\end{equation*}
Consider the following adjoint equations:
\begin{equation*}
\left\{
    \begin{array}{l}
    dr_{\rho}^{\varepsilon}(t) = \left[ c_y(t)r_{\rho}^{\varepsilon}(t) - b_y(t)p_{\rho}^{\varepsilon}(t) - \sigma_y(t)q_{\rho}^{\varepsilon}(t) - \sum_{i=0}^l \lambda_{\rho}^{i,\varepsilon} f_y^i(t) \right] dt\\
    \hspace{3cm} \vspace{0.2cm}+ \left[ c_z(t)r_{\rho}^{\varepsilon}(t) - b_z(t)p_{\rho}^{\varepsilon}(t) - \sigma_z(t)q_{\rho}^{\varepsilon}(t) - \sum_{i=0}^l \lambda_{\rho}^{i,\varepsilon} f_z^i(t) \right] dW_t\\
    \vspace{0.2cm} dp_{\rho}^{\varepsilon}(t) = - \left[ c_x(t)r_{\rho}^{\varepsilon}(t) - b_x(t)p_{\rho}^{\varepsilon}(t) - \sigma_x(t)q_{\rho}^{\varepsilon}(t) - \sum_{i=0}^l \lambda_{\rho}^{i,\varepsilon} f_x^i(t) \right] dt + q_{\rho}^{\varepsilon}(t) dW_t\\
    \vspace{0.2cm} r_{\rho}^{\varepsilon}(0) = -\mathbb{E} \left[ \sum_{i=0}^l \lambda_{\rho}^{i,\varepsilon} g^i(y(0)) \right], \hspace{0.2cm} p_{\rho}^{\varepsilon}(T) = -\varphi_x(x(T))r_{\rho}^{\varepsilon}(T) + \sum_{i=0}^l \lambda_{\rho}^{i,\varepsilon} h_x^i(x(T))
    \end{array}
\right.
\end{equation*}
By the duality relation in \cite{MR2230926} (in the proof of Theorem 4.1, p. 168) we then get that
\begin{equation*}
    \mathbb{E} \left[ \int_0^T \delta \mathcal{H}_{\rho} \chi_{E_{\varepsilon}}(t) dt \right] \leq \sqrt{\rho} \varepsilon + o(\varepsilon)
    \label{IntroIntHamEst}
\end{equation*}
where the hamiltonian $\mathcal{H}_{\rho}$ is
\begin{equation*}
    \mathcal{H}_{\rho}(t,x,y,z,u,r,p,q, \lambda_{\rho}^{0,\varepsilon}, \lambda_{\rho}^{\varepsilon}) = p \cdot b - r \cdot c + q \cdot \sigma - \sum_{i=0}^l\lambda_{\rho}^{i,\varepsilon} f^i,
\end{equation*}
and further, by (\ref{IntroFactor1}) and (\ref{IntroFactor2}) and a property of Clarke's generalized gradient $\partial d_{\Lambda}$ (see Lemma \ref{AppendixMetric2Lemma}), there exists a choice of multipliers $(\lambda_{\rho}^{0,\varepsilon}, \lambda_{\rho}^{1,\varepsilon},\ldots, \lambda_{\rho}^{l,\varepsilon}  )$ such that
\begin{equation*}
    \lambda_{\rho}^{0,\varepsilon} \geq o(1), \hspace{0.5cm} |\lambda_{\rho}^{0,\varepsilon}|^2 + |\lambda_{\rho}^{\varepsilon}|^2 = 1 + o(1).
\end{equation*}
Thus, there is a subsequence, still denoted by $(\lambda _{\rho}^{0,\varepsilon}, \lambda_{\rho}^{\varepsilon})$, such that
\begin{equation*}
    \smash{\displaystyle\lim_{\varepsilon \rightarrow 0}} (\lambda _{\rho}^{0,\varepsilon}, \lambda_{\rho}^{\varepsilon}) = (\lambda _{\rho}^{0}, \lambda_{\rho})
\end{equation*}
for some $(\lambda _{\rho}^{0}, \lambda_{\rho}) \in \R^{1+l}$ such that
\begin{equation*}
    \lambda_{\rho}^0 \geq 0, \hspace{0.5cm} |\lambda_{\rho}^0|^2 + |\lambda_{\rho}|^2 = 1.
    \label{IntroLagMultRho}
\end{equation*}
Furthermore, by \cite{MR1613669}, we have strong convergence in $L_{\mathcal{F}}^2(\Omega; C([0,T]; \mathbb{R})) \times L_{\mathcal{F}}^2(\Omega; C([0,T]; \mathbb{R})) \times L_{\mathcal{F}}^2(0,T;\mathbb{R})$ of (adapted) solutions to the FBSDEs:
\begin{equation*}
    (r_{\rho}^{\epsilon}(\cdot), p_{\rho}^{\epsilon}(\cdot), q_{\rho}^{\epsilon}(\cdot)) \rightarrow (r_{\rho}(\cdot), p_{\rho}(\cdot), q_{\rho}(\cdot)).
\end{equation*}
For the limit properties of the hamiltonian, let $E_{\varepsilon} = [\bar{t}, \bar{t} + \varepsilon]$ for an arbitrary $\bar{t} \in [0,T]$ and let $A$ be an arbitrary set in $\mathcal{F}_{\bar{t}}$ and consider the control $u(t) = v \mathbb{1}_A + u_{\rho}(t) \mathbb{1}_{\Omega \backslash A}$,
for an arbitrary $v \in U$ in (\ref{IntroIntHamEst}). Dividing (\ref{IntroIntHamEst}) by $\varepsilon$ and then approaching the limit $\varepsilon \rightarrow 0$ yields
\begin{equation*}
    \begin{array}{l}
    \sqrt{\rho} \geq \mathcal{H}_{\rho}(t,x(t),y(t),z(t),u(t),r(t),p(t),q(t), \lambda_{\rho}^{0,\varepsilon}, \lambda_{\rho}^{\varepsilon}) \vspace{0.2cm}\\
    \hspace{1.3cm}- \mathcal{H}_{\rho}(t,x(t),y(t),z(t),u_{\rho}(t),r(t),p(t),q(t), \lambda_{\rho}^{0,\varepsilon}, \lambda_{\rho}^{\varepsilon})
    \end{array}
\end{equation*}
Finally, let $\rho \rightarrow 0$. By choosing a subsequence and using (\ref{IntroLagMultRho}) we may assume that $(\lambda_{\rho}^0, \lambda_{\rho}) \rightarrow (\lambda^0, \lambda) \in \R^{1+l}$
for which (\ref{IntroLagMult}) holds. Furthermore, by the fact that $\bar{d}(u_{\rho}(\cdot), \bar{u}(\cdot)) \leq \sqrt{\rho}$ and the theorem on continuous dependence of a parameter in FBSDEs in \cite{MR1613669} we have strong convergence of $(r_{\rho}(\cdot), p_{\rho}(\cdot), q_{\rho}(\cdot))$ to a triple $(r(\cdot), p(\cdot), q(\cdot))$ satisfying (\ref{IntroFBSDEConstSMP}). The transversality condition (\ref{IntroTransversalityCond}) follows in the limit $\varepsilon , \rho \rightarrow 0$ from the definition of $\lambda_{\rho}^{\varepsilon}$ and $\partial d_{\Lambda}(\cdot)$ since then we have the inequality
    $\langle \lambda_{\rho}^{\varepsilon}, v - \mathbb{E} [\mathbf{x}_{\rho}^{\varepsilon}(T)] \rangle \geq - o(1), \hspace{0.2cm} \text{ as } \varepsilon \rightarrow 0.$
\end{proof}
\begin{remark}
As in Remark \ref{IntroRemarkThmSMP1}, analogue principles hold in Theorem \ref{IntroConstrainedSMPthm}.
\end{remark}

\section{The Problem}
\label{theproblem}
The present literature of Principal-Agent problems in continuous time covers several different (although related) models. Basically two of these are fundamental and of particular interest; The \textit{Full Information} case and the \textit{Hidden Action} case. In this section we develop a method of characterizing optimal contracts in the Hidden Action regime, based solely on the framework presented in Section \ref{Preliminaries}.\\
\indent The common setup of the Principal-Agent problem under Hidden Action is inspired by the seminal paper of Holmstr\"{o}m-Milgrom (see \cite{MR882097}) and is well treated, for instance, in \cite{MR2963805} and \cite{Williams}. Because of the considerable difficulty arising in such models the usual approach is to state the problem in the \textit{weak form}, which is a mathematically related (but not equivalent) problem. The sole reason for this mathematical difficulty lies at the heart of the problem formulation, namely the asymmetric flow of information between the principal and the agent. \\
\indent In what follows we first consider a model under Hidden Action deviating from classical formulation by allowing the Principal to observe more information than what is usually done (the full filtration generated by the noise). The advantage of this being that we get a problem which is mathematically tractable in the strong formulation, whereas the disadvantage is the immediate reduction of realism. We will refer to this model as the \textit{Hidden Contract} model. After this we explain how the weak formulation of the classical Principal-Agent problem is contained in the Hidden Contract framework. We end the section with a discussion of the models and future prospects.\\

\subsection{Hidden Contract in the strong formulation }
\label{HiddenContract}

\indent Consider a Principal-Agent model where output $x(t)$ is modelled as a risky asset solving the SDE:
\begin{equation}
\left\{
    \begin{array}{l}
    dx(t) = f(t, x(t), e(t)) dt + \sigma (t, x(t)) dW_t, \hspace{0.5cm} t\in (0,T], \\
    x(0) = 0,
    \end{array}
\right.
\end{equation}
Here $T > 0$ and $W_t$ is a 1-dimensional standard Brownian motion defined on the filtered probability space $(\Omega, \mathcal{F}, \mathbb{F}, \mathbb{P})$. The functions $f$ and $\sigma$ represent production rate and volatility respectively, and we assume both of them satisfy the regularity conditions (S2)-(S3) from the previous section. The process $e(\cdot)$ represents the agent's level of effort, taking values in some predefined subset $E \subseteq \R$ (typically $E = [0,\hat{e}]$ for some non-negative $\hat{e}$, or $E = \R$) and is required to belong to $\mathcal{E}[0,T]$, where
\begin{equation*}
    \mathcal{E}[0,T] := \{ e: [0,T] \times \Omega \rightarrow E; \hspace{0.2cm} e \text{ is } \mathcal{F}_t \text{-adapted} \}
\end{equation*}
We consider the case of Hidden Actions meaning that the principal cannot observe the effort $e(\cdot)$. Output, however, is always public information and observed by both the principal and the agent. In this model we relax the public information structure to include also the Brownian noise, i.e. the filtration $\mathbb{F}$ rather than $\mathbb{F}_x \subset \mathbb{F}$. Thus, apart from the effort policy, the principal and the agent have access to the same information. Before the period starts the principal specifies an $\mathbb{F}$-adapted cash-flow $s(\cdot)$ (typically non-negative) for all $t \in [0,T]$, thus compensating the agent for costly effort in managing $x(\cdot)$ and by that providing incentives. Just as for the effort we assume $s(t) \in S$ for all $t \in [0,T]$ and some subset $S \subseteq \R$ and require $s(\cdot) \in \mathcal{S}[0,T]$, where
\begin{equation*}
    \mathcal{S}[0,T] := \{ s: [0,T] \times \Omega \rightarrow S;  \hspace{0.2cm} s \text{ is } \mathcal{F}_t \text{-adapted} \}.
\end{equation*}
The principal is not constrained by any means and can, in principle, commit to any process $s(\cdot) \in \mathcal{S}[0,T]$.\\
\indent In this model we consider cost functionals $\mathcal{J}_P$ and $\mathcal{J}_A$ of the principal and the agent respectively of the following form:
\begin{equation}
    \mathcal{J}_A (e(\cdot); s) := \E \left[ \int_0^T u (t,x(t),e(t),s(t)) dt + v(x(t)) \right]
    \label{TheProblemAgentCost}
\end{equation}
and
\begin{equation}
    \mathcal{J}_P(s(\cdot)) := \E \left[ \int_0^T \mathcal{U} (t,x(t),s(t)) dt + \mathcal{V}(x(t)) \right].
\end{equation}
The agent will accept $s(\cdot)$ and start working for the principal only if it fulfills the participation constraint:
\begin{equation}
    \mathcal{J}_A (\bar{e}(\cdot );s) \leq W_0,
    \label{TheProblemPartConst}
\end{equation}
for some, typically negative, constant $W_0$. By assuming incentive compatibility, meaning that the agent will act as to optimize $\mathcal{J}_A$ in response to any given $s(\cdot)$, we may think of the Principal-Agent problem as divided into two coupled problems; \textit{The Agent's problem} and \textit{The Principal's problem}.\\
\noindent \textbf{The Agent's problem:}
Given any $s(\cdot) \in \mathcal{S}[0,T]$ (which we assume do satisfy the participation constraint (\ref{TheProblemPartConst})) the Agent's problem is to find a process $\bar{e}(\cdot) \in \mathcal{E}[0,T]$ such that
\begin{equation*}
    \mathcal{J}_A (\bar{e}(\cdot );s) = \smash{\displaystyle\inf_{e(\cdot) \in \mathcal{E}[0,T]}} \mathcal{J}_A(e(\cdot);s).
\end{equation*}\\

\noindent \textbf{The Principal's problem:}
Given that the Agent's problem has an optimal solution $\bar{e}(\cdot)$ the Principal's problem is to find a process $\bar{s}(\cdot) \in \mathcal{S}[0,T]$ such that
\begin{equation*}
    \mathcal{J}_P (\bar{s}(\cdot)) = \smash{\displaystyle\inf_{s(\cdot) \in \mathcal{S}[0,T]}} \mathcal{J}_P(s(\cdot))
\end{equation*}
and
\begin{equation*}
    \mathcal{J}_A (\bar{e}(\cdot );\bar{s}) \leq W_0.
\end{equation*}\\
\noindent In this context the following definition is natural.
\begin{definition}
An \textit{optimal contract} is a pair $( \bar{e}(\cdot), \bar{s}(\cdot) ) \in \mathcal{E}[0,T] \times \mathcal{S}[0,T]$ obtained by sequentially solving the Agent's- and the Principal's problem.
\end{definition}
\noindent In the language of game theory an optimal contract can thus be thought of as a Stackelberg-equilibrium in a two-player non-zero-sum game.\\
\begin{remark}
In the classical formulation of the problem the principal is only allowed to write a contract based on his/her observation of output. Thus, rather than being adapted to $\mathbb{F}$ the process $s(\cdot)$ must be adapted to the smaller filtration $\mathbb{F}_x$. In that case the most general structure of a cash-flow is $s_t=s(t,(x)_t)$, where $(x)_t$ denotes the history of output in the time interval $[0,t]$. The restriction of having $s_t$ adapted to $\mathbb{F}_x$ complicates the mathematical treatment in the strong formulation and one has to consult the framework of weak solutions of SDEs. In the above formulation of the model we do allow for a general cash-flow to be $\mathbb{F}$-adapted. Since the agent cannot control the noise he/she has no use of the precise structure of $s(\cdot)$ and we may therefore regard it as hidden information, thus motivating the name \textit{Hidden Contract}. The agent influence the contract through \ref{TheProblemPartConst} and reacts to it by choosing optimal effort as the project starts.\\
\indent For applications it is of course not realistic to assume a reasonable contract to be purely $\mathbb{F}$-adapted. However, as we develop the scheme of solving the Hidden Contract model it becomes clear that many cases allow for $\mathbb{F}_x$-adapted optimal cash-flows $s(\cdot)$. This is in line with present literature in Full Information-models (e.g. \cite{MR2237178}) and motivates our approach.
\end{remark}
\indent It is important to note that even though the principal cannot observe the agent's effort, he/she can still offer the agent a contract by suggesting a choice of effort $e(\cdot)$ and a compensation $s(\cdot)$. However, by incentive compatibility the principal knows that the agent will only follow such a contract if the suggested effort solves the agent's problem. In order to find the optimal suggested effort, $\bar{e}(\cdot)$, the principal must therefore have information of the agent's preferences, i.e. the functions $u$ and $v$. The realism of such an assumption is indeed questionable but nevertheless necessary in our formulation due to the participation constraint. In order to make the intuition clear and to avoid any confusion we adopt the convention that the principal has full information of the agent's preferences $u$ and $v$. Of course, this does not serve any mathematical purpose, it solely gives a tractable way of thinking of how actual contracting is realized.\\
\indent Thus, the principal can predict any optimal effort $\bar{e}(\cdot)$ of the agent's problem and by that suggest an optimal contract $(\bar{e}(\cdot), \bar{s}(\cdot))$, if it exists.\\
\indent The idea is now to apply the methods from Section \ref{Preliminaries} to characterize optimal contracts in the general Principal-Agent model presented above. The agent's problem is standard in stochastic optimal control of SDEs, with Hamiltonian
\begin{equation}
    \mathcal{H}_A (t, x, e, p, q, s) := p \cdot f(t, x, e) + q \cdot \sigma(t, x) - u(t, x, e, s).
    \label{TheProblemAgentHam}
\end{equation}
Therefore, by Theorem \ref{IntroThmSMP1} we have for any optimal pair $(\bar{x}(\cdot), \bar{e}(\cdot))$ of the Agent's problem the existence of a pair of adjoint processes $(p(\cdot), q(\cdot))$ solving the BSDE:
\begin{equation}
\left\{
    \begin{array}{l}
        dp(t) = - \left\{ f_x(t, \bar{x}(t), \bar{e}(t)) p(t) + \sigma_x(t, \bar{x}(t)) q(t) - u_x (t, \bar{x}(t), \bar{e}(t)) \right\} dt + q(t) dW_t, \\
        p(T) = - v_x(\bar{x}(T)),
    \end{array}
    \label{TheProblemAgentBSDE}
\right.
\end{equation}
for which
\begin{equation}
    \bar{e}(t)= \smash{\displaystyle\arg\max_{e \in E}}\hspace{0.1cm}\mathcal{H}_A(t, \bar{x}(t), e, p(t), q(t), s(t)),
    \label{TheProblemAgentHamMax}
\end{equation}
for a.e. $t \in [0,T]$ and $\Prob$-a.s.

\noindent In many applications relation (\ref{TheProblemAgentHamMax}) is sufficient in order to find a closed form expression for $\bar{e}(t)$ directly, but not always. Before proceeding to the Principal's problem we assume such a closed form expression and write
\begin{equation*}
    \bar{e}_t = e^*(t, \bar{x}(t), p(t), q(t), s(t)),
\end{equation*}
where $e^* : \R_+ \times \R^4 \rightarrow \R$ is a function having sufficient regularity to allow for existence of a unique solution to the FBSDE (\ref{TheProblemFBSDE}) below.
\begin{remark}
In the above formulation the agent is reacting to a cash-flow $s(\cdot)$ given solely as an $\mathcal{F}_t$-adapted process $s: [0,T] \times \Omega \rightarrow S$. Thus the more refined structure of $s$, for instance its dependence of $x$, is unknown to the agent. The case in which $s$ is given as a function of output is briefly discussed in Section~\ref{PaymentByResultContracting}.
\end{remark}
\begin{remark}
The precise regularity theory of $e^*$ following from the general setup of the Agent's problem is rather involved. In particular, consider the special case when: $v(x) \equiv 0$, $u \in C^3$, $u_{ee}(\cdot) > 0$, $u(\cdot, e)/|e| \geq \gamma(e)$ where $\gamma(e) \rightarrow \infty$ as $|e| \rightarrow \infty$, $u(\cdot) \geq -c_1$ and $u(\cdot)|_{e=0} \leq c_1$ for a suitable constant $c_1$, $|u_x(\cdot)| \leq c_2 u(\cdot) + c_3$ for suitable constants $c_2$ and $c_3$, and $|u_e(\cdot)| \leq C(R)$ whenever $|e| \leq R$ for suitable $C(R)$. Then $e^*$ is a \textit{Markov control policy}, i.e. a continuous function in all variables, being Lipschitz continuous on discs and of linear growth in $x$ (see e.g. \cite{MR1199811}, Theorem 11.1 p. 206).
\end{remark}
Based on the information given by $e^*$ the principal wishes to minimize the cost $\mathcal{J}_P$ by selecting a process $s(\cdot)$ respecting (\ref{TheProblemPartConst}). The dynamics of the corresponding control problem is, in contrast to the simple SDE of the agent's problem, a FBSDE built up by the output SDE coupled to the agent's adjoint BSDE. More precisely:
\begin{equation}
    \left\{
    \begin{array}{l}
        d\bar{x}_t = f(t, \bar{x}(t), e^*(t, \bar{x}(t), p(t), q(t), s(t))) dt + \sigma (t, \bar{x}(t)) dW_t,\\
        dp(t) = - \left\{ f_x(t, \bar{x}(t), e^*(t, \bar{x}(t), p(t), q(t), s(t))) p(t) + \sigma_x(t, \bar{x}(t)) q(t) \right.\\
        \hspace{4cm}\left. - u_x (t, \bar{x}(t), e^*(t, \bar{x}(t), p(t), q(t), s(t))) \right\} dt + q(t) dW_t,\\
        \bar{x}(0) = 0, \hspace{0.2cm} p(T) = - v_x(\bar{x}(T)).
    \end{array}
    \right.
    \label{TheProblemFBSDE}
\end{equation}
Since the participation constraint (\ref{TheProblemPartConst}) is closed and convex (i.e. the constraint set $\{ x \in \R : x \leq W_0 \}$ is closed and convex) we can, under suitable regularity of the coefficients in (\ref{TheProblemFBSDE}), apply Theorem \ref{IntroConstrainedSMPthm} to characterize optimality in the Principal's problem. Therefore, with the Hamiltonian of the principal given by
\begin{equation}
    \begin{array}{l}
    \mathcal{H}_P (t,x,p,q,s,R,P,Q,\lambda_P,\lambda_A) :=\\
    \hspace{0.5cm}-R \cdot \left\{ f_x(t, x, e^*(t, x, p, q, s)) p + \sigma_x(t, x) q - u_x (t, x, e^*(t, x, p, q, s)) \right\}\\
    \hspace{0.5cm}+ P \cdot f(t, x, e^*(t, x, p, q, s)) + Q \cdot \sigma (t, x) - \lambda_A u (t, x, e^*(t, x, p, q, s))\\
    \hspace{0.5cm}- \lambda_P  \mathcal{U} (t,x,s),
    \end{array}
    \label{TheProblemPrincipalHam}
\end{equation}
we have for any optimal 4-tuple $(\bar{x}(\cdot), \bar{p}(\cdot),\bar{q}(\cdot), \bar{s}(\cdot))$ of the Principal's problem the existence of Lagrange multipliers $\lambda_A, \lambda_P \in \R$ satisfying
\begin{equation*}
    \lambda_P \geq 0, \hspace{0.5cm} \lambda_A^2 + \lambda_P^2 = 1,
    \label{HiddenContracLagrange}
\end{equation*}
and a triple of adjoint processes $(R(\cdot), P(\cdot), Q(\cdot))$ solving the FBSDE:
\begin{equation}
\left\{
        \begin{array}{l}
        \vspace{0.2cm}dR_t = \left\{ \left[\partial_p f_{x}(t, \bar{x}(t), e^*(t, \bar{x}(t), \bar{p}(t), \bar{q}(t), \bar{s}(t))) p(t) \right.\right.\\
        \vspace{0.2cm}\hspace{2.8cm}\left.+ f_{x}(t, \bar{x}(t), e^*(t, \bar{x}(t), \bar{p}(t), \bar{q}(t), \bar{s}(t)))\right.\\
        \vspace{0.2cm}\hspace{2.8cm}\left.- \partial_{p} u_x (t, \bar{x}(t), e^*(t, \bar{x}(t), \bar{p}(t), \bar{q}(t), \bar{s}(t)),\bar{s}(t)) \right] \cdot R_t\\
        \vspace{0.2cm}\hspace{1cm}\left.- \partial_p f (t, \bar{x}(t), e^*(t, \bar{x}(t), \bar{p}(t), \bar{q}(t), \bar{s}(t)))\cdot P(t)\right.\\
        \vspace{0.2cm}\hspace{2.8cm}\left.+ \lambda_A \partial_p u (t, \bar{x}(t), e^*(t, \bar{x}(t), \bar{p}(t), \bar{q}(t), \bar{s}(t)),\bar{s}(t))\right\} dt\\
        \vspace{0.2cm}\hspace{0.7cm}+\left\{ \left[ \partial_q f_x(t, \bar{x}(t), e^*(t, \bar{x}(t), \bar{p}(t), \bar{q}(t), \bar{s}(t))) \bar{p}(t) + \sigma_x(t, \bar{x}(t))\right.\right.\\
        \vspace{0.2cm}\hspace{2.8cm}\left.\left.- \partial_{q} u_x (t, \bar{x}(t), e^*(t, \bar{x}(t), \bar{p}(t), \bar{q}(t), \bar{s}(t)),\bar{s}(t)) \right] \cdot R_t \right.\\
        \hspace{1cm}\vspace{0.2cm}\left.- \partial_q f(t, \bar{x}(t), e^*(t, \bar{x}(t), \bar{p}(t), \bar{q}(t), \bar{s}(t))) \cdot P(t)\right.\\
        \vspace{0.4cm}\hspace{2.8cm}\left.+ \lambda_A \partial_q u (t, \bar{x}(t), e^*(t, \bar{x}(t), \bar{p}(t), \bar{q}(t), \bar{s}(t)),\bar{s}(t)) \right\} dW_t,\\
        \vspace{0.2cm}dP(t) = -\left\{ - \left[ \partial_x f_x(t, \bar{x}(t), e^*(t, \bar{x}(t), \bar{p}(t), \bar{q}(t), \bar{s}(t))) \bar{p}(t) + \sigma_{xx}(t, \bar{x}(t)) \bar{q}(t)\right.\right.\\
        \hspace{2.8cm}\vspace{0.2cm}\left.\left.- \partial_{x} u_x (t, \bar{x}(t), e^*(t, \bar{x}(t), \bar{p}(t), \bar{q}(t), \bar{s}(t)),\bar{s}(t)) \right] \cdot R_t\right.\\
        \vspace{0.2cm}\hspace{1cm}\left.+ \partial_x f(t, \bar{x}(t), e^*(t, \bar{x}(t), \bar{p}(t), \bar{q}(t), \bar{s}(t)) \cdot P(t)+ \sigma_x (t, \bar{x}(t)) \cdot Q(t)\right.\\
        \hspace{1cm}\vspace{0.2cm}\left. - \lambda_A \partial_{x} u(t, \bar{x}(t), e^*(t, \bar{x}(t), \bar{p}(t), \bar{q}(t), \bar{s}(t)),\bar{s}(t)) - \lambda_P \mathcal{U}_x (t,\bar{x}(t),s(t)) \right\} dt\\
        \hspace{1cm}\vspace{0.2cm}+ Q(t) dW_t,\\
        R(0) = 0,\hspace{0.2cm} P(T) =  v_{xx}(\bar{x}(T))R(T) - \lambda_A v_x (\bar{x}(T))- \lambda_P \mathcal{V}_x(\bar{x}(T)),
        \end{array}
        \label{TheProblemPrincipalFBSDE}
    \right.
\end{equation}
so that
\begin{equation*}
        \bar{s}(t)= \smash{\displaystyle\arg\max_{s \in S}}\hspace{0.1cm} \mathcal{H}_P(t, \bar{x}(t),\bar{p}(t),\bar{q}(t),s,R(t),P(t),Q(t),\lambda_P,\lambda_A),
\end{equation*}
$\text{a.e. } t \in [0,T], \hspace{0.3cm} \Prob\text{-a.s.}$ The notation $\partial_z$ in (\ref{TheProblemPrincipalFBSDE}) means the total derivative with respect to variable $z$.
\noindent We are now ready to state a full characterization of optimal contracts in the Hidden Contract Principal-Agent problem. Before doing so we introduce the following technical assumption:
\begin{description}
    \item[(PA)] All functions involved in the Agent's problem satisfy the conditions (S1)-(S3) from Section \ref{Preliminaries}. The functions defining the Principal's problem (including composition with the map $e^*$) satisfy the conditions (SC1)-(SC3), also from Section \ref{Preliminaries}.
\end{description}
\begin{theorem}[Hidden Contract Optimality Characterization]
    Let the regularity conditions in (PA) hold and consider a Principal-Agent problem under Hidden Contracts in continuous time with participation constraint defined by $W_0 < 0$. Then, if $(\bar{e}(\cdot), \bar{s}(\cdot))$ is an optimal contract there exist numbers $\lambda_A, \lambda_P \in \R$ such that
    \begin{equation*}
        \lambda_P \geq 0, \hspace{0.5cm} \lambda_A^2 + \lambda_P^2 = 1,
    \end{equation*}
    a pair $(p(\cdot), q(\cdot))\in L_{\mathcal{F}}^2(0,T;\R) \times (L_{\mathcal{F}}^2(0,T;\R))$ solving the SDE in (\ref{TheProblemAgentBSDE}) and a triple $(R(\cdot), P(\cdot), Q(\cdot)) \in L_{\mathcal{F}}^2(\Omega; C([0,T]; \mathbb{R})) \times L_{\mathcal{F}}^2(\Omega; C([0,T]; \mathbb{R})) \times L_{\mathcal{F}}^2(0,T;\mathbb{R})$ solving the FBSDE in (\ref{TheProblemPrincipalFBSDE}) such that, sequentially,
\begin{equation*}
    \bar{e}(t) = \smash{\displaystyle\arg\max_{e \in E}}\hspace{0.1cm}\mathcal{H}_A(t, \bar{x}(t), e, p(t), q(t), s(t)),
\end{equation*}
and
\begin{equation*}
    \bar{s}(t)= \smash{\displaystyle\arg\max_{s \in S}}\hspace{0.1cm} \mathcal{H}_P(t, \bar{x}(t),\bar{p}(t),\bar{q}(t),s,R(t),P(t),Q(t),\lambda_P,\lambda_A),
\end{equation*}
    with Hamiltonians $\mathcal{H}_A$ and $\mathcal{H}_P$ as in (\ref{TheProblemAgentHam}) and (\ref{TheProblemPrincipalHam}) respectively.
    \label{TheProblemMainThm}
\end{theorem}
\begin{remark}
    The regularity conditions in (PA) may indeed seem too restrictive for most applications. It is, however, important to remember that the purpose of (PA) is to ensure unique solvability of all stochastic differential equations involved in the general setting. An alternative way of formulating Theorem \ref{TheProblemMainThm} is therefore to exchange (PA) by simply assuming the existence and uniqueness of solutions to (\ref{TheProblemAgentBSDE}) and (\ref{TheProblemPrincipalFBSDE}).
\end{remark}

\subsection{Hidden Action in the weak formulation}
\label{HiddenAction}
Turning to the classical theory of Principal-Agent problems under Hidden Action we now assume the information available to the principal to be strictly contained in the filtration generated by output $\mathbb{F}_x$. The set of admissible cash-flows therefore has to be redefined in comparison with Section \ref{HiddenContract} to the subset $\mathcal{S}_x[0,T] \subset \mathcal{S}[0,T]$, where
\begin{equation}
    \mathcal{S}_x[0,T] := \{ s: [0,T] \times \Omega \rightarrow S;  \hspace{0.2cm} s \text{ is } \mathcal{F}_t^x \text{-adapted} \}.
\end{equation}
A direct approach (in the spirit of Section \ref{HiddenContract} for the Hidden Contract model) is not mathematically tractable in the Agent's problem, and so far no author has been able to carry this through. Therefore, in line with \cite{MR2963805} and \cite{Williams} we make the problem tractable via weak solutions of SDEs. First we model the output as
\begin{equation}
\left\{
    \begin{array}{l}
        dx(t) = \sigma(t,x(t)) dW_t,\\
        x(0) = 0,
    \end{array}
\right.
\end{equation}
where we assume that $\sigma > 0$ and $\mathbb{E}\left[ \int_0^t \sigma(t,x(t))^2 dt \right] < \infty$. We then consider the process $\Gamma^e(t)$ which solves the SDE
\begin{equation}
\left\{
    \begin{array}{l}
        d\Gamma^e(t) = \Gamma^e(t) f(t, x(t), e(t)) \sigma^{-1}(t, x(t)) dW_t,\\
        \Gamma^e(0) = 0,
    \end{array}
    \label{HiddenActionAgentDyn}
\right.
\end{equation}
where $e(\cdot)$ is the effort of the agent and $f(t,x,e)$ is as in the previous section. Note that $\Gamma^e = \mathcal{E}(Y)$ where $Y(t) = \int_0^t f(s,x(s),e(s))\cdot \sigma^{-1}(s,x(s))dW_s$ and $\mathcal{E}(\cdot)$ denotes the Dolean's stochastic exponential. The key idea behind the weak formulation of the Hidden Action model is to let the agent control the \textit{density of output} $\Gamma^e(t)$ rather than the output process itself, and to consider $(x)_T$ as a fixed but random realization (actually $\mathbb{F}_x = \mathbb{F}$ as a consequence of the regularity of $\sigma$). Indeed, assuming e.g. Novikov's condition for $\Gamma^e$ we have by Girsanov's theorem that the probability measure $d\mathbb{P}^e$ defined by
\begin{equation}
    \frac{d\mathbb{P}^e}{d\mathbb{P}} = \Gamma^e(T)
\end{equation}
makes the process $W^e(t)$ defined by
\begin{equation}
    dW^e_t = dW_t - f(t,x(t),e(t))\sigma^{-1}(t,x(t))dt
\end{equation}
a $\mathbb{P}^e$-Brownian motion. In particular
\begin{equation}
    dx(t) = f(t,x(t),e(t))dt + \sigma(t,x(t))dW^e_t
    \label{HiddenActionWeakDyn}
\end{equation}
and
\begin{equation}
    \begin{array}{l}
        \vspace{0.2cm}\mathcal{J}_A(e(\cdot);s) = \mathbb{E}^e \left[ \int_0^T u(t,x(t),e(t),s(t))dt v(x(T)) \right] =\\
        \hspace{0.3cm}= \mathbb{E} \left[ \int_0^T \Gamma^e(t)u(t,x(t),e(t),s(t))dt + \Gamma^e(T)v(x(T)) \right].
    \end{array}
    \label{HiddenActionAgentCost}
\end{equation}
In the weak formulation of the Hidden Action model the Agent's- and the Principal's problems are formulated as follows.
\vspace{0.2cm}\\
\noindent \textbf{The Agent's problem (weak formulation):}
Given any $s(\cdot) \in \mathcal{S}_x[0,T]$ (fulfilling the participation constraint) the Agent's problem is to find a process $\bar{e}(\cdot) \in \mathcal{E}[0,T]$ such that cost functional is
\begin{equation*}
    \mathcal{J}_A (\bar{e}(\cdot );s) = \mathbb{E} \left[ \int_0^T \Gamma^e(t)u(t,x(t),e(t),s(t))dt + \Gamma^e(T)v(x(T)) \right],
\end{equation*}
is minimized, subject to the dynamics in (\ref{HiddenActionAgentDyn}).\vspace{0.2cm}\\
\noindent \textbf{The Principal's problem (strong formulation):}
Given that the Agent's problem has an optimal solution $\bar{e}(\cdot)$ in the weak formulation the Principal's problem is to find a process $\bar{s}(\cdot) \in \mathcal{S}[0,T]$,
such that the cost functional
\begin{equation*}
    \mathcal{J}_P (\bar{s}(\cdot)) := \mathcal{J}_P(s(\cdot)) := \E \left[ \int_0^T \mathcal{U} (t,x(t),s(t)) dt + \mathcal{V}(x(t)) \right]
\end{equation*}
is minimized and
\begin{equation*}
    \mathcal{J}_A (\bar{e}(\cdot );\bar{s}) = \mathbb{E} \left[ \int_0^T u(t,x(t),e(t),s(t))dt + v(x(T)) \right] \leq W_0.
\end{equation*}
subject to the dynamics
\begin{equation*}
\left\{
    \begin{array}{l}
    dx(t) = \sigma (t, x(t)) dW_t, \hspace{0.5cm} t\in (0,T], \\
    x(0) = 0.
    \end{array}
\right.
\end{equation*}
\begin{remark}
    Here we have chosen to formulate the Principal's problem in the strong form rather than in the weak form, which seems to be most common in the literature. However, as pointed out in \cite{MR2963805}, due to adaptivity this approach can be problematic in certain models. A fact that one should be aware of.
\end{remark}
With the full Principal-Agent problem stated as above we can apply the same technique for solving the Hidden Action case as we did with the Hidden Contract case. However, since the control variable $e$ figures in the diffusion coefficient of (\ref{HiddenActionAgentDyn}) we require the following convexity assumption in order to avoid the second order adjoint processes in the maximum principle:\vspace{0.2cm}\\
\noindent \textbf{Assumtion 1:} The set $E \subset \R$ is convex.\vspace{0.2cm}\\
The Agent's Hamiltonian in the weak formulation is
\begin{equation}
    \mathcal{H}_A(t,x,\Gamma^e,e,p,q,s) := q \cdot \Gamma^e \cdot \frac{f(t, x, e)}{\sigma(t, x)} - \Gamma^e \cdot u(t,x,e,s),
    \label{HiddenActionAgHam}
\end{equation}
and by Theorem \ref{IntroThmSMP1} any optimal control $\bar{e}(t)$ solving the Agent's problem must maximise $\mathcal{H}_A$ where the pair $(p(\cdot), q(\cdot))$ solves the Agent's adjoint BSDE:
\begin{equation}
\left\{
    \begin{array}{l}
        dp(t) =  -\left\{ q(t) \cdot \frac{f(t, x(t), \bar{e}(t))}{\sigma(t, x(t))} - u(t,x(t),\bar{e}(t),s(t)) \right\}dt + q(t) dW_t,\\
        p(T) = -v(x(T))
    \end{array}
    \right.
    \label{HiddenActionAdjBSDE}
\end{equation}
If the functions $f$ and $u$ both are differentiable in the $e$ variable and we assume that $\bar{e}(\cdot) \in \text{int}(E)$, maximizing $\mathcal{H}_A$ translates into the first order condition
\begin{equation}
    q(t) = \sigma(t,x(t)) \cdot \frac{u_e(t,x(t),\bar{e}(t), s(t))}{f_e(t,x(t),\bar{e}(t))},
\end{equation}
which is in agreement with \cite{Williams}.\\
\indent Just as in the Hidden Contract case we now proceed into the Principal's problem by assuming the existence of a closed form expression for optimal effort $\bar{e}_t = e^*(t, x(t), q(t), s(t))$ (having sufficient regularity to allow for existence and uniqueness of a solution to (\ref{HiddenActionFBSDE})). The Principal is facing the problem of minimizing the cost functional $\mathcal{J}_P$ by controlling the following FBSDE:
\begin{equation}
\left\{
    \begin{array}{l}
        \vspace{0.2cm}dx(t) = \sigma (t, x(t)) dW_t,\\
        \vspace{0.2cm}dp(t) =  - \left\{ q(t) \cdot \frac{f(t, x(t), e^*(t, x(t), q(t), s(t)))}{\sigma(t, x(t))} - u(t,x(t),e^*(t, x(t), q(t), s(t)),s(t)) \right\} dt\\
        \vspace{0.2cm}\hspace{1.3cm}+ q(t) dW_t,\\
        x(0) = 0, p(T) = -v(x(T)).
        \label{HiddenActionFBSDE}
    \end{array}
\right.
\end{equation}
We now apply Theorem \ref{IntroConstrainedSMPthm} in order to characterize optimal cash-flows in the Principal's problem. The associated Hamiltonian reads
\begin{equation}
    \begin{array}{l}
    \vspace{0.2cm}\mathcal{H}_P(t,x,q,s,R,P,Q,\lambda_P,\lambda_A) := R \cdot \left\{ -q \cdot \frac{f(t, x, e^*(t, x, q, s))}{\sigma(t, x)} + u(t,x,e^*(t, x, q, s),s) \right\}\\
    \hspace{0.1cm} + Q \cdot \sigma (t, x) - \lambda_A \cdot u(t,x,e^*(t, x, q, s),s) - \lambda_P \cdot \mathcal{U} (t,x,s),
    \end{array}
    \label{HiddenActionPrincHam}
\end{equation}
and for any optimal 4-tuple $(\bar{x}(\cdot), \bar{p}(\cdot), \bar{q}(\cdot), \bar{s}(\cdot))$ we have the existence of Lagrange multipliers $\lambda_A, \lambda_P \in \R$ satisfying (\ref{HiddenContracLagrange}) and a triple of processes $(R(\cdot), P(\cdot), Q(\cdot))$ solving the adjoint FBSDE:
\begin{equation}
\left\{
    \begin{array}{l}
        \vspace{0.2cm}dR(t) = \displaystyle\left\{ \left[\frac{f(t,x(t),e^*(t,x(t),q(t),s(t))) + q(t) \cdot \partial_q f(t,x(t),e^*(t,x(t),q(t),s(t)))}{\sigma(t,x(t))} \right.  \right.\\
        \vspace{0.2cm}\hspace{4cm} + \partial_q u(t,x(t),e^*(t,x(t),q(t),s(t)),s(t)) \bigg] \cdot R(t)\\
        \vspace{0.3cm}\hspace{4.9cm}+ \lambda_A \partial_q u(t,x(t),e^*(t,x(t),q(t),s(t)),s(t)) \bigg\} dW_t\\

        \vspace{0.2cm}dP(t) = \displaystyle-\left\{ \left[ q(t) \cdot \frac{\partial_x f(t,x(t),e^*(t,x(t),q(t),s(t))) \cdot \sigma(t,x(t))}{\sigma(t,x(t))^2}  \right. \right.\\
        \vspace{0.2cm}\hspace{1.9cm}  \displaystyle - q(t) \cdot \frac{ f(t,x(t),e^*(t,x(t),q(t),s(t))) \cdot \sigma_x(t,x(t))}{\sigma(t,x(t))^2}\\
        \vspace{0.2cm}\hspace{1.9cm}+ \partial_x u(t,x(t),e^*(t,x(t),q(t),s(t)),s(t)) \bigg] \cdot R(t) + \sigma_x(t,x(t)) \cdot Q(t)\\
        \vspace{0.2cm}\hspace{1.7cm}-\lambda_A \partial_x u(t,x(t),e^*(t,x(t),q(t),s(t)),s(t)) - \lambda_P \mathcal{U}_x(t,x(t),s(t)) \bigg\} dt\\
        \vspace{0.3cm}\hspace{10.2cm}+Q(t) dW_t\\

        R(0) = 0, \hspace{0.2cm}P(T) = (R(T) - \lambda_A) \cdot v_x(x(T)) - \lambda_P V_x(x(T))
    \end{array}
    \right.
    \label{HiddenActionAdjFBSDE}
\end{equation}
We are now ready to state the full characterization of optimal contracts in the Hidden Action case in line with Theorem \ref{TheProblemMainThm}.
\begin{theorem}
    Let the regularity conditions in (PA) hold and consider a Principal-Agent problem under Hidden Actions in continuous time with participation constraint defined by $W_0 < 0$. Then, if $(\bar{e}(\cdot), \bar{s}(\cdot))$ is an optimal contract there exist numbers $\lambda_A, \lambda_P \in \R$ such that
        \begin{equation*}
            \lambda_P \geq 0, \hspace{0.5cm} \lambda_A^2 + \lambda_P^2 = 1,
        \end{equation*}
    a pair $(p(\cdot), q(\cdot))\in L_{\mathcal{F}}^2(0,T;\R) \times (L_{\mathcal{F}}^2(0,T;\R))$ solving the SDE in (\ref{HiddenActionAdjBSDE}) and a triple $(R(\cdot), P(\cdot), Q(\cdot)) \in L_{\mathcal{F}}^2(\Omega; C([0,T]; \mathbb{R})) \times L_{\mathcal{F}}^2(\Omega; C([0,T]; \mathbb{R})) \times L_{\mathcal{F}}^2(0,T;\mathbb{R})$ solving the FBSDE in (\ref{HiddenActionAdjFBSDE}) such that, sequentially,
        \begin{equation*}
            \bar{e}(t) = \smash{\displaystyle\arg\max_{e \in E}}\hspace{0.1cm}\mathcal{H}_A(t,\bar{x}(t),\Gamma^e(t),e,p(t),q(t),s(t)),
        \end{equation*}
    and
        \begin{equation*}
            \bar{s}(t)= \smash{\displaystyle\arg\max_{s \in S}}\hspace{0.1cm} \mathcal{H}_P(t,\bar{x}(t),\bar{q}(t),s,R(t),P(t),Q(t),\lambda_P,\lambda_A),
        \end{equation*}
    with Hamiltonians $\mathcal{H}_A$ and $\mathcal{H}_P$ as in (\ref{HiddenActionAgHam}) and (\ref{HiddenActionPrincHam}) respectively.
    \label{HiddenActionMainThm}
\end{theorem}

\subsection{Conclusion}
In this section we have studied two different models of Principal-Agent type. The first model, referred to as the Hidden Contract model, is to the best of our knowledge new whereas the second model of Hidden Action in the weak formulation is by now well known. We believe that our contribution to the existing literature is of interest for two reasons:
\begin{itemize}
  \item The Hidden Contract model is indeed deviating from the main stream direction of research in Principal-Agent problems. One of the most serious claims against our model is probably the lack of realism since the Agent has to accept a hidden contract based on a cash-flow that does not necessarily depend of the observed output. Since the Agent can only influence the contract via the participation constraint he/she may well be exposed to undesirable high levels of risk without being aware of it. The realism in our Hidden Contract model would therefore increase notably if the participation constraint could be extended to include also a statement of acceptable risk, for instance by imposing that the variance of the accumulated cash-flow is bounded from above. This will be the topic of a follow-up paper as it requires extending  the approach of Section \ref{theproblem} and the framework in Section \ref{Preliminaries} to cost-functionals of mean-field type.\vspace{0.3cm}

  \item The Hidden Action case in the weak formulation as stated above does not add any significant contribution to the existing theory. As mentioned in the Introduction our approach is closest to the one in \cite{Williams} since both use stochastic generalizations of Pontryagin's maximum principle to characterize optimality. The main difference compared to \cite{Williams} is that we exploit these methods to find a full characterization covering both the Agent's and the Principal's problem, and not only that of the Agent.\\
      \indent The book \cite{MR2963805} presents a thorough characterization of the Hidden Action model, adopting the convention of keeping the weak formulation also for the Principal's problem. The solution scheme is then carried through via a well known technique in stochastic optimization which is slightly different from our approach. However, since our next step is to consider Principal-Agent problems involving nonlinear cost functionals such as variance we believe that the method described in Section \ref{HiddenAction} is more suitable.
\end{itemize}

\section{A Solved Example}
\begin{figure}[h!]
  \centering
\includegraphics[width=1.0\linewidth]{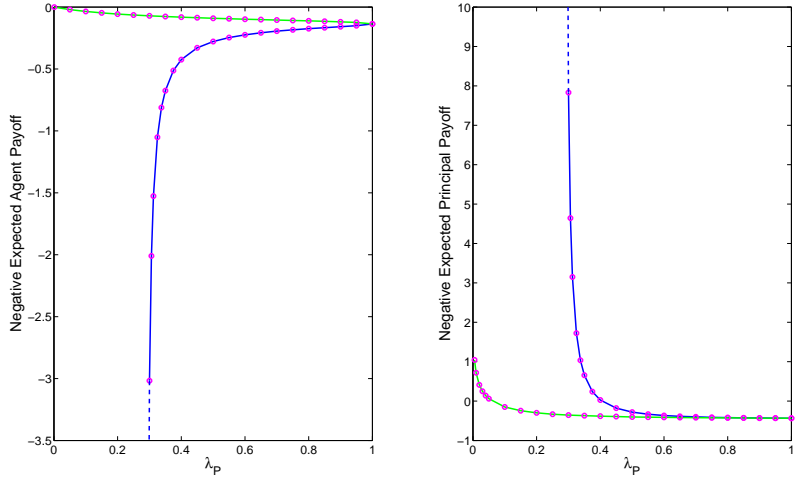}%
  \caption{Monte-Carlo simulations of $\mathcal{J}_A(\bar{e}(\cdot);\bar{s})$ and $\mathcal{J}_P(\bar{s}(\cdot))$ as functions of $\lambda_P$ based on $5 \cdot 10^5$ sample paths at each point. Parameter values: $a = b = \sigma = 1, \alpha = 0.2, \beta = 1, T = 0.35$. Blue signifies the branch $\lambda_A > 0$ whereas green is the branch $\lambda_A < 0$.}
  \label{LQexpectedpayofffig}
\end{figure}

To illustrate the method described in Section \ref{theproblem} we devote this section to solve a concrete example of the Hidden Contract Principal-Agent problem in a linear-quadratic (LQ)-setting. It turns out that closed form solutions of both the optimal effort $\bar{e}(t)$ and the optimal cash-flow $\bar{s}(t)$ can be found. Consider the following dynamics of production,
\begin{equation*}
    \left\{
    \begin{array}{l}
        dx(t) = (ax(t) + be(t))dt + \sigma dW_t, \hspace{0.3cm} t\in (0,T],\\
        x(0) = 0, \hspace{0.3cm} a,b \in \R \text{ and } \sigma > 0,
    \end{array}
    \right.
\end{equation*}
and let the preferences of the agent and the principal be described by quadratic utility functions as:
\begin{eqnarray}
    && \mathcal{J}_A (e(\cdot);s) := \E\left[ \int_0^T \frac{(s_t - e_t)^2}{2} dt - \alpha \cdot \frac{x(T)^2}{2} \right],\\
    && \mathcal{J}_P (s(\cdot)) := \E\left[\int_0^T \frac{s_t^2}{2} dt - \beta \cdot \frac{x(T)^2}{2}  \right].
    \label{LQpayoffP}
\end{eqnarray}
Note that we are following the convention of Section \ref{theproblem} and consider cost- rather than payoff-functionals. Thus, the Agent's utility function should be interpreted as a desire to maintain a level of effort close to the compensation given by the cash-flow.
We think of the parameters $\alpha>0$ and $\beta>0$ as bonus factors of total production at time $T$. For the participation constraint we require any admissible cash-flow $s(t)$ to satisfy
\begin{equation}
    \mathcal{J}_A(\bar{e}(\cdot); s) \leq W_0, \hspace{0.5cm} W_0 < 0,
    \label{LQconstr}
\end{equation}
where $\bar{e}(\cdot)$ denotes the optimal response of the agent given $s(t)$.\\
\indent Assume that the principal offers the agent $s(t)$ over the period $0 \leq t \leq T$. The Hamiltonian function of the agent's problem is
\begin{equation*}
    \mathcal{H}_A(x,e,p,q,s) := p \cdot (ax + be) + q \cdot \sigma - \frac{(s-e)^2}{2},
\end{equation*}
so
\begin{equation}
    \frac{\partial \mathcal{H}_A}{\partial e} = bp + s - e = 0 \hspace{0.5cm} \text{and} \hspace{0.5cm} \bar{e}(t)= bp(t) + s(t),
    \label{LQopteff}
\end{equation}
where the pair $(p,q)$ solves the adjoint equation
\begin{equation*}
    \left\{
    \begin{array}{l}
        dp(t) = -ap(t) dt + q(t) dW_t,\\
        p(T) = \alpha x(T).
    \end{array}
    \right.
\end{equation*}
\noindent Turning now to the principal's problem we would like to control the FBSDE
\begin{equation}
    \left\{
    \begin{array}{l}
        dx(t) = (ax(t) - b^2p(t) + bs(t))dt + \sigma dW_t,\\
        dp(t) = -ap(t)dt + q(t)dW_t,\\
        x(0)=0, p(T)=\alpha x(T).
    \end{array}
    \label{LQFBSDE1}
    \right.
\end{equation}
optimally by $s$ with respect to (\ref{LQpayoffP}) under the state constraint in (\ref{LQconstr}) that now, given the optimal effort $\bar{e}$ in (\ref{LQopteff}), reads $\mathcal{J}_A(bp(\cdot) + s(\cdot)) \leq W_0$. Thus the Hamiltonian function of this problem is given by
\begin{equation*}
    \begin{array}{l}
        \mathcal{H}_P (x,p,q,s,P,Q,R,\lambda_A,\lambda_P) :=\\
         \vspace{0.2cm}\hspace{0.7cm}= P \cdot (ax + b^2p + bs) - R \cdot ap + Q \cdot \sigma - \lambda_A \cdot \frac{b^2p^2}{2} - \lambda_P \cdot \frac{s^2}{2},
    \end{array}
\end{equation*}
where $\lambda_A$ and $\lambda_P$ are Lagrange-multipliers such that $\lambda_P \geq 0$ and $\lambda_P^2 + \lambda_A^2 = 1$, according to Theorem \ref{IntroConstrainedSMPthm} in Section \ref{Preliminaries}. Again we easily find the unique candidate for optimal control:
\begin{equation*}
    \frac{\partial \mathcal{H}_P}{\partial s} = bP - \lambda_P s, \hspace{0.5cm} \text{so} \hspace{0.5cm} \bar{s}(t)= \frac{b}{\lambda_P}P(t)
\end{equation*}
where the triple $(R(t),P(t),Q(t))$ solves the adjoint FBSDE:
\begin{equation}
    \left\{
    \begin{array}{l}
        dR(t) = (aR(t) - b^2P(t) + \lambda_a b^2 p(t))dt,\\
        dP(t) = -aP(t) dt + Q(t) dW_t,\\
        R(0) = 0, P(T) = -\alpha R(T) + (\alpha \lambda_A + \beta \lambda_P)x(T).
    \end{array}
    \right.
    \label{LQFBSDE2}
\end{equation}
In order to find concrete expressions for the controls $(\bar{e}_t,\bar{s}_t)$ in the optimal contract we must solve the BSDEs in (\ref{LQFBSDE1}) and (\ref{LQFBSDE2}). For this we consider a general ansatz of feedback type
\begin{equation*}
    p(t) = \varphi (t,x(t),R(t)), \hspace{0.5cm} \text{ and } \hspace{0.5cm} P(t) = \psi (t,x(t),R(t))
\end{equation*}
and compute the stochastic differentials by utilizing It\^{o}'s lemma. We get
\begin{equation*}
    \begin{array}{l}
        dp(t) = d\varphi (t,x(t),R(t)) =\\
        \vspace{0.2cm}\left\{ \partial_t \varphi + (ax(t) + b^2\varphi + \frac{b^2}{\lambda_P} \psi) \varphi_x + (aR(t)+ \lambda_A b^2 \varphi - b^2 \psi) \varphi_y + \frac{\sigma ^2}{2} \varphi_{xx} \right\} + \sigma \varphi_xdW_t,
    \end{array}
\end{equation*}
\begin{equation*}
    \begin{array}{l}
        dP(t) = d\psi (t,x(t),R_t) =\\
        \vspace{0.2cm}\left\{ \partial_t \psi + (ax(t) + b^2\varphi + \frac{b^2}{\lambda_P} \psi) \psi_x + (aR_t + \lambda_A b^2 \varphi - b^2 \psi) \psi_y + \frac{\sigma ^2}{2} \psi_{xx} \right\} dt + \sigma \psi_xdW_t.
    \end{array}
\end{equation*}
By identifying coefficients in the equations we end up with the following system of semilinear parabolic PDE's:
\begin{equation}
\left\{
    \begin{array}{l}
        \partial_t \varphi + a\varphi + (ax(t) + b^2\varphi + \frac{b^2}{\lambda_P} \psi) \varphi_x + (aR_t + \lambda_A b^2 \varphi - b^2 \psi) \varphi_y + \frac{\sigma ^2}{2} \varphi_{xx} = 0,\\
        \partial_t \psi + a\psi + (ax(t) + b^2\varphi + \frac{b^2}{\lambda_P} \psi) \psi_x + (aR_t + \lambda_A b^2 \varphi - b^2 \psi) \psi_y + \frac{\sigma ^2}{2} \psi_{xx} = 0,\\
        \varphi(T,x,y) = \alpha x, \psi(T,x,y) = -\alpha y + (\alpha \lambda_a + \beta \lambda_P)x,\\
        \vspace{0.2cm}a,b,\lambda_A,\lambda_P \in \R, \hspace{0.1cm} \alpha \in (0,1), \hspace{0.1cm} \sigma > 0, \hspace{0.1cm} \lambda_P \geq 0, \hspace{0.1cm} \lambda_P^2 + \lambda_A^2 = 1.
    \end{array}
\right.
\label{LQsemilinparsys}
\end{equation}
This system is of course in general very difficult or even impossible to solve explicitly, but a numerical scheme like for instance finite differences would handle the problem well. However, in this case we may proceed analytically a bit further by considering yet another ansatz:
\begin{equation*}
    \varphi(t,x,y) = C_1(t) x + C_2(t) y, \hspace{0.5cm} \psi(t,x,y) = D_1(t) x + D_2(t) y.
\end{equation*}
By identifying coefficients, this ansatz reduces (\ref{LQsemilinparsys}) to a terminal conditioned Riccati-system:
\begin{equation}
\left\{
    \begin{array}{l}
        \frac{dC_1}{dt} =  2aC_1(t) + b^2C_1^2(t) + \frac{b^2}{\lambda_P}C_1(t)D_1(t) - b^2C_2(t)D_1(t) + \lambda_A b^2 C_1(t)C_2(t),\\
        \frac{dC_2}{dt} =  2aC_2(t) + \lambda_A b^2 C_2^2(t) - b^2 C_2(t)D_2(t) + \frac{b^2}{\lambda_P} C_1(t)D_2(t) + b^2 C_1(t)C_2(t),\\
        \frac{dD_1}{dt} =  2aD_1(t) + \frac{b^2}{\lambda_P} D_1^2(t) + b^2 C_1(t)D_1(t) + \lambda_A b^2 C_1(t)D_2(t) - b^2 D_1(t)D_2(t),\\
        \frac{dD_2}{dt} =  2aD_2(t) - b^2D_2^2(t) + \lambda_A b^2 C_2(t)D_2(t) + b^2D_1(t)C_2(t) + \frac{b^2}{\lambda_P} D_1(t)D_2(t),\\
        \vspace{0.2cm} C_1(T) = \alpha, \hspace{0.1cm} C_2(T) = 0, \hspace{0.1cm} D_1(T) = \alpha \lambda_A + \beta \lambda_P, \hspace{0.1cm} D_2(T) = -\alpha.
    \end{array}
\right.
\label{LQriccatisyst}
\end{equation}
Riccati-systems have closed form solutions only in exceptional cases and the one in (\ref{LQriccatisyst}) is too general to admit such (even though the coefficients are constant). A good theoretical treatment of Riccati-systems and matrix Riccati equations can be found in \cite{MR1997753}.\\
\indent Even though an explicit solution seems hopeless, it is an easy task to solve (\ref{LQriccatisyst}) numerically, given the parameters. An example is presented in Fig. \ref{LQriccatisystfig} below.
\begin{figure}[h!]
  \centering
\includegraphics[width=0.8\linewidth]{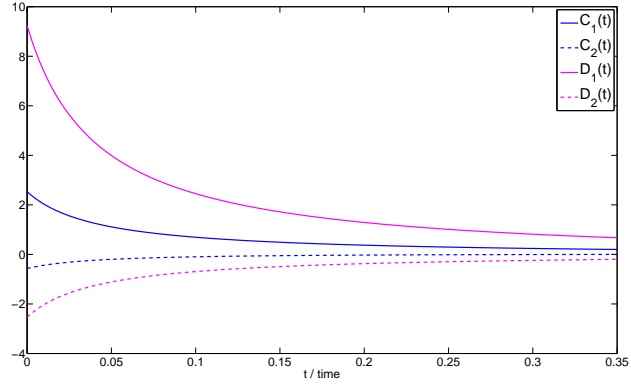}%
  \caption{Solution curves of (\ref{LQriccatisyst}) with parameter values chosen as: $a = b = \sigma = 1, \alpha = 0.2, \beta = 1, \lambda_P = 0.5, \lambda_A = +\sqrt{1-\lambda_p^2}, T = 0.35.$}
  \label{LQriccatisystfig}
\end{figure}
Hence, for the LQ Principal-Agent we get the unique semi-explicit optimal contract as
\begin{equation}
        \begin{array}{l}
        \{\bar{e}(t), \bar{s}(t)\} =\vspace{0.2cm}\\
        \hspace{0.4cm}= \{ (D_1(t) + bC_1(t)) \bar{x}(t) + (D_2(t) + bC_2(t))\bar{R}(t), \frac{b}{\lambda_P}D_1(t) \bar{x}_t + \frac{b}{\lambda_P}D_2(t)\bar{R}(t) \}
        \end{array}
        \label{LQOptCon}
\end{equation}
where the optimal dynamics $(\bar{x}(t), \bar{R}(t))$ solve the forward SDE's in (\ref{LQFBSDE1}) and (\ref{LQFBSDE2}), replacing $p(t)$ and $P(t)$ by $\varphi(t,\bar{x}(t),\bar{R}(t))$ and $\psi(t, \bar{x(t)}, \bar{R}(t))$. What still remains is to find a feasible pair $(\lambda_A, \lambda_P)$ so that the optimal contract fulfills the participation constraint in (\ref{LQconstr}). One way of finding such a pair is for instance by stochastic simulation of $(\bar{x}(t), \bar{R}(t))$ (e.g. a simple Euler-Maruyama scheme) and then estimate the payoff in (\ref{LQconstr}) by Monte-Carlo techniques for different values of $\lambda_P$. In Fig. \ref{LQexpectedpayofffig} below we have included the results of such a scheme where the negative expected payoff for the agent and the principal are plotted as functions of the multiplier $\lambda_P$.
Note that $\bar{R}(t)$ satisfies the linear ODE
\begin{equation}
\left\{
\begin{array}{l}
    \frac{dR}{dt} + (b^2 D_2(t) - \lambda_A b^2 C_2(t) - a) R(t) = (\lambda_A b^2 C_1(t) - b^2 D_1(t)) \bar{x}(t)\\
    R(0) = 0,
    \end{array}
    \right.
\end{equation}
so
\begin{equation*}
    R(t) = \frac{\int_0^t \text{exp} \left\{ \int_0^s b^2 D_2(u) - \lambda_A b^2 C_2(u) - a \hspace{0.1cm} du \right\} \cdot (\lambda_A b^2 C_1(s) - b^2 D_1(s)) \bar{x}(t) ds }{\text{exp} \left\{ \int_0^t b^2 D_2(s) - \lambda_A b^2 C_2(s) - a \hspace{0.1cm} ds \right\}},
\end{equation*} 
and is by that $\mathbb{F}_x$-adapted. Therefore, in this model the optimal contract (\ref{LQOptCon}) is $\mathbb{F}_x$-adapted and coincides with the corresponding strong solution to the Hidden Action problem, i.e. when the information set of the Principal is generated by output.\\
\indent An interesting observation to make from the numerical result is the monotonic nature of the agents negative payoff. Apparently $\mathcal{J}_A(\bar{e}(\cdot); \bar{s}))$ is a decreasing function of $\lambda_A \in (-1,1)$ ranging over the entire half-axis $(-\infty,0)$, meaning that the agent increases his/her expected final utility by an increased value of $\lambda_A$. In contrast, the function $\mathcal{J}_P(\bar{s}(\cdot))$ is a non-monotonic function of $\lambda_A \in (-1,1)$ but does increase on each of the branches $-1 < \lambda_A \leq 0$ and $0 \leq \lambda_A < 1$ (considering $\lambda_A = 0$ as the starting point). Furthermore (as suggested by the numerics) the negative branch is strictly below the positive one. Thus, for any $W_0 \in (W_c, 0]$ where $W_c := \mathcal{J}_A(\bar{e}(\cdot); \bar{s})|_{\lambda_A = 0}$ and corresponding multiplier $\lambda_A^0$ such that $\mathcal{J}_A(\bar{e}(\cdot);\bar{s})|_{\lambda_A = \lambda_A^0} = W_0$ we have that
\begin{equation*}
\left\{
    \begin{array}{l}
    \mathcal{J}_A(\bar{e}(\cdot);\bar{s})|_{\lambda_A = 0} < W_0,\\
    \mathcal{J}_P(\bar{s}(\cdot))|_{\lambda_A = 0} < \mathcal{J}_P(\bar{s}(\cdot))|_{\lambda_A = \lambda_A^0},
    \end{array}
\right.
\end{equation*}
meaning that both parts will benefit from a contract offering the agent higher expected utility than the minimal threshold $-W_0$. This implication is somewhat counterintuitive since it would be reasonable to expect the optimal contract to fulfill the equality
\begin{equation}
    \mathcal{J}_A(\bar{e}(\cdot);\bar{s}) = W_0.
    \label{LQcounterintuitive}
\end{equation}
Indeed, in the case of lump-sum paying contracts with separable utility functions the relation (\ref{LQcounterintuitive}) is proven in \cite{MR2465709} in a continuous time setting. The effects suggested above therefore reveal an interesting additional complexity in optimal contracting with non-separable utility functions and continuously paying cash-flows. To the best of our knowledge this has not been observed in the present literature. It has in fact rather conjectured false in the setting considered in \cite{MR2465709}.

\section{A View Towards Hidden Actions in the Strong Formulation}
\label{PaymentByResultContracting}
So far in our presentation of the Principal-Agent problem under Hidden Action we have been working under the crucial assumption that the
Agent is facing a stochastic optimal control problem in the weak sense. However, from the point of view of applications it is probably more natural to consider the strong formulation in which the the agents controls output rather than the density of output. This problem is much more difficult and we devote this section to a case study, only to give a glimpse of the type of difficulties that may arise. We also present a possible direction for future research in this area.\\
\indent Just as in Section \ref{theproblem} we consider a hidden action Principal-Agent model where the output $x(t)$ is a risky asset solving the SDE
\begin{equation}
\left\{
    \begin{array}{l}
    dx(t) = f(t, x(t), e(t)) dt + \sigma (t, x(t)) dW_t,\hspace{0.5cm} t\in (0,T],\\
    x(0) = 0,
    \end{array}
\right.
\end{equation}
and the cost functionals are of the form
\begin{equation}
    \mathcal{J}_A (e(\cdot);s) := \E \left[ \int_0^T u (t,x(t),e(t),s(t, x(t))) dt + v(x(t)) \right],
\end{equation}
and
\begin{equation}
    \mathcal{J}_P(s(\cdot)) := \E \left[ \int_0^T \mathcal{U} (t,x(t),s(t)) dt + \mathcal{V}(x(t)) \right].
\end{equation}
The principal wishes to find an optimal contract $(\bar{e}(\cdot), \bar{s}(\cdot))$ solving the Principal-Agent problem by tailoring a differentiable (in $x$) map $(t, x) \mapsto s(t, x)$, i.e. a 'smooth' cash-flow of feed-back type. Compared to Section \ref{theproblem} essential differences now arise already in the Agent's problem. For instance the agent's BSDE takes the form
\begin{equation*}
\left\{
    \begin{array}{l}
        dp(t) = - \left\{ f_x(t, \bar{x}(t), \bar{e}(t)) p(t) + \sigma_x(t, \bar{x}(t)) q(t) - \partial_x u (t, \bar{x}(t), \bar{e}(t), s(t, x(t))) \right\} dt \\ \qquad\quad + q(t) dW_t,\\
        p(T) = - v_x(\bar{x}(T)),
    \end{array}
\right.
\end{equation*}
and the Hamiltonian is
\begin{equation*}
    \mathcal{H}_A (t, x, e, p, q, s) := p \cdot f(t, x, e) + q \cdot \sigma(t, x) - u(t, x, e, s(t,x)).
\end{equation*}
The extra dependence on $x(t)$ appearing in these equations may seem harmless, but it actually complicates the problem substantially. Having the additional information of how the cash-flow $s(\cdot)$ depends on the output makes it possible for the agent to select optimal effort more carefully by estimating future dividends of the contract. Loosely speaking this yields a path dependence of $s(\cdot)$ in $\bar{e}(\cdot)$, turning the Principal's problem into a non-standard stochastic optimal control problem. In order to illustrate this we consider an explicit example.
\begin{example}
Consider the dynamics governing output as
\begin{equation*}
    \left\{
    \begin{array}{l}
        dx(t) = e(t) dt + \sigma dW_t,\\
        x(0) = 0
    \end{array}
    \right.
\end{equation*}
and cost functionals
\begin{equation}
    \mathcal{J}_A (e(\cdot);s) := \E \left[ \int_0^T \frac{e(t)^2}{2} - s(t,x(t))dt \right]
\end{equation}
and
\begin{equation}
    \mathcal{J}_P(s(\cdot)) := \E \left[ \int_0^T s(t)dt - x(T) \right].
    \label{PaymentByResultContractingPrincCost}
\end{equation}
We let the participation constraint be as in (\ref{LQconstr}) and consider the Agent's problem. Clearly
\begin{equation*}
    \mathcal{H}_A (t, x, e, p, q, s) := p \cdot e + q \cdot \sigma - \left( \frac{e^2}{2} - s(t,x) \right)
\end{equation*}
so
\begin{equation*}
    \frac{\partial \mathcal{H}_A}{\partial e} = p - e = 0 \hspace{0.2cm} \text{and} \hspace{0.2cm} \bar{e}_t = p(t)
\end{equation*}
where $(p(\cdot), q(\cdot))$ solves the BSDE
\begin{equation}
\left\{
    \begin{array}{l}
        dp(t) = - s_x(t, x(t)) dt + q(t) dW_t,\\
        p(T) = 0.
    \end{array}
\right.
\label{PaymentByResultContractingLQBSDE}
\end{equation}
In order to solve (\ref{PaymentByResultContractingLQBSDE}) we consider the ansatz $p(t) = \varphi (t, x(t))$, for a smooth function $\varphi$, and apply It\^{o}'s lemma so that
\begin{equation*}
    dp(t) = \left\{ \partial_t \varphi(t,x(t)) + \varphi(t,x(t)) \varphi_x(t,x(t)) + \frac{\sigma^2}{2} \varphi_{xx} (t,x(t))\right\} dt + \sigma \varphi_x (t,x(t))dW_t
\end{equation*}
Thus, by identifying coefficients we get the equation
\begin{equation}
\left\{
    \begin{array}{l}
        \partial_t \varphi(t,x) + \varphi(t,x) \varphi_x(t,x) + \frac{\sigma^2}{2} \varphi_{xx} (t,x)+ s_x(t,x) = 0,\\
        \varphi(T, x) = 0,
    \end{array}
\right.
\label{PaymentByResultContractingBurger}
\end{equation}
i.e. Burger's equation with force term $s_x$. Despite the non-linear nature of this PDE it admits an explicit solution via the so called Hopf-Cole transformation
\begin{equation*}
    \varphi (t,x) = \frac{\partial}{\partial x} \left[ \sigma^2 \text{log} v(t,x) \right],
\end{equation*}
where $v(t,x)$ solves the inhomogeneous heat equation
\begin{equation}
\left\{
    \begin{array}{l}
        \partial_t v(t,x)= -\frac{\sigma^2}{2} v_{xx}(t,x) - s(t,x),\\
        v(T,x) = \text{exp}\left( \frac{1}{\sigma^2} \int \varphi(T,x) dx \right) = 1.
    \end{array}
    \label{PaymentByResultContractingTerValProb}
\right.
\end{equation}
By time-reversal we turn (\ref{PaymentByResultContractingTerValProb}) into an initial value problem and thereby find the solution to (\ref{PaymentByResultContractingBurger}) as
\begin{equation}\label{phi}
    \varphi (t,x,s_.) = \sigma^2 \frac{\int_t^T \int_{\R} \frac{\xi - x}{\sigma^2(\tau - t)} G(x, \xi, \tau - t) s(\tau, \xi) d\xi d\tau}{1 + \int_t^T \int_{\R} G(x, \xi, \tau - t) s(\tau, \xi) d\xi d\tau},
\end{equation}
where $G(x, \xi, \tau - t)$ is the Gauss-kernel
\begin{equation*}
    G(x, \xi, \tau - t) := \frac{1}{\sqrt{2\pi \sigma^2 (\tau - t)}}\text{exp} \left\{ -\frac{(\xi - x)^2}{2\sigma^2 (\tau - t)}\right\}.
\end{equation*}
We use the notation '$s_.$' in (\ref{phi}) to underline the {\it functional dependence} of the function $s(\cdot)$ in the solution $\varphi$.\\

\noindent It is easy to check that if $s(t,x) = C$ for some constant $C > 0$, then $\bar{e}_t~=~\varphi(t, x(t))~\equiv~0$, which of course is to expect. Another interesting example is to take $s(t,x) = x$, making the optimal effort non-trivial:
\begin{equation*}
    \bar{e}(t) = \frac{\sigma^2 (T-t)}{1 + x(t)(T-t)}.
\end{equation*}
Having $\bar{e}(t) = \varphi (t, x(t))$ the Principal's problem is to control
\begin{equation}
    \left\{
    \begin{array}{l}
        dx(t) = \varphi(t,x(t), s_.) dt + \sigma dW_t,\\
        x(0) = 0,
    \end{array}
    \right.
    \label{PaymentByResultContractingPrincDyn}
\end{equation}
optimally by minimizing the cost functional (\ref{PaymentByResultContractingPrincCost}), yet respecting the participation constraint
\begin{equation*}
    \mathcal{J}_A (\varphi(\cdot); s_.) = \E \left[ \int_0^T \frac{\varphi(t,x(t),s_.)^2}{2} - s(t,x(t))dt \right] < W_0,
\end{equation*}
for some given $W_0 < 0$. The main obstacle as compared to the general theory in Section \ref{theproblem} is the dependence of $s(\cdot)$ in (\ref{PaymentByResultContractingPrincDyn}) now being functional rather than point-wise. To handle this one could for instance consider the probability density function $\mu(t, z)$ of $x(t)$ at time $t \geq 0$. Since (\ref{PaymentByResultContractingPrincDyn}) is an It\^{o}-process with initial condition $x(0) = 0$ we know that $\mu(t, z)$ solves the Fokker-Planck equation (or, the Kolmogorov forward equation):
\begin{equation}
\left\{
    \begin{array}{l}
        \frac{\partial}{\partial t} \mu(t, z) = - \frac{\partial}{\partial z} \left[ \varphi(t, z, s_.) \mu(t, z) \right] + \frac{\sigma^2}{2}\frac{\partial^2}{\partial z^2} \mu(t, z)\\
        \mu(t, z) = \delta_0
    \end{array}
    \label{PaymentByResultContractingFokPla}
\right.
\end{equation}
Having $\mu(t, z)$ at hand the Principal's problem can now be formulated as a deterministic problem in optimal control of PDEs. More precisely, the principal wishes to find a function $s(t,x)$ (lying in a suitable function space making (\ref{PaymentByResultContractingFokPla}) well-posed) minimizing the functional
\begin{equation*}
    \mathcal{J}_P(s(\cdot)) = \int_0^T \int_{\R} s(t, z) \mu(t, z) dz dt - \int_{\R} z \mu(T, z) dz
\end{equation*}
subject to the constraint
\begin{equation*}
    \mathcal{J}_A (\varphi(\cdot); s_.) = \int_0^T \int_{\R} \left( \frac{\varphi(t,z, s_.)^2}{2} - s(t,z) \right) \mu(t, z) dz dt < W_0.
\end{equation*}
This problem, however, goes beyond the scope of our presentation and we leave it for a future study. In \cite{MR2724518} and recently in \cite{Annunziato} the authors explore the general possibility of the above approach to stochastic optimal control, thus formulating the problem in terms of optimal control of the Fokker-Planck equation. Their setting is, however, significantly different from ours in that it does not allow for a functional dependence of the control in (\ref{PaymentByResultContractingFokPla}).
\end{example}

\appendix
\section{Technical Theorems}
\label{AppendixTechThms}
The main result upon which the proof of Theorem \ref{IntroConstrainedSMPthm} relies is the following well-known theorem by Ekeland (see e.g. \cite{MR1696772}, Lemma 6.2 p. 145).
\begin{theorem}[Ekeland's variational principle]
Let $(V,d)$ be a complete metric space and let $F : V \rightarrow ( -\infty , \infty ]$ be a proper (i.e. $F \equiv +\infty$), lower semicontinuous function bounded from below. Let $v_0 \in \mathcal{D} := \{ v \in V | F(v) < \infty \}$ and let $\lambda > 0$ be fixed. Then there exists a $\bar{v} \in V$ such that:
\begin{description}
    \item[(i)] $F(\bar{v}) + \lambda d(\bar{v}, v_0) \leq F(v_0)$\vspace{0.2cm}
    \item[(ii)] $F(\bar{v}) < F(v) + \lambda d(v, \bar{v}), \hspace{0.5cm} \forall v \neq \bar{v}.$
\end{description}
\label{AppendixEkelandThm}
\end{theorem}
\noindent The next result is a useful consequence of Theorem \ref{AppendixEkelandThm}.
\begin{corollary}
Let the assumptions of Theorem \ref{AppendixEkelandThm} hold. Let $\rho > 0$ and $v_0 \in V$ be such that
\begin{equation*}
    F(v_0) \leq \smash{\displaystyle\inf_{v \in V}} F(v) + \rho.
\end{equation*}
Then there exists a $v_{\rho} \in V$ such that
\begin{equation*}
    F(v_{\rho}) \leq F(v_0), \hspace{0.5cm} d(v_{\rho}, v_0) \leq \sqrt{\rho},
\end{equation*}
and for all $v \in V$,
\begin{equation*}
    -\sqrt{\rho} d(v,v_{\rho}) \leq F(v) - F(v_{\rho}).
\end{equation*}
\label{AppendixEkelandCor}
\end{corollary}
\noindent In addition to Ekeland's variational principle we also need a couple of lemmas to prove Theorem \ref{IntroConstrainedSMPthm}. One dealing with a certain metric structure on the space $\mathcal{U}[0,T]$, and the other explaining some regularity properties of the canonical distance function $d_{\Lambda}(\cdot, \cdot)$ to a given closed and convex set $\Lambda$. Both results are presented with proofs in \cite{MR1696772} (Lemma 6.4 p. 147 and Lemma 6.5 p. 148).
\begin{lemma}
Let $\bar{d}$ be defined by the following:
\begin{equation*}
\bar{d}(u(\cdot), \tilde{u}(\cdot)) := |\{ (t,\omega) \in [0,T] \times \Omega: u(t,\omega) \neq \tilde{u}(t,\omega) \}|, \hspace{0.3cm} \forall u(\cdot), \tilde{u}(\cdot) \in \mathcal{U}[0,T],
\end{equation*}
where $|A|$ denotes the product measure of the Lebesgue measure and the probability $\mathbb{P}$ of a set $A \subseteq [0,T] \times \Omega$. Then $\bar{d}$ is a metric under which $\mathcal{U}[0,T]$ is a complete metric space.
\label{AppendixMetricLemma}
\end{lemma}
\begin{lemma}
Let $\Lambda \subseteq \R^l$ be a closed and convex set and define the distance function
\begin{equation*}
    d_{\Lambda} (v) := \smash{\displaystyle\inf_{v' \in \Lambda}} |v - v'|, \hspace{0.5cm} \forall v \in \R^l.
\end{equation*}
Then the following holds:
\begin{equation*}
\begin{array}{l}
    \textbf{(i) } d_{\Lambda} : \R^l \rightarrow \R \textit{ is convex and Lipschitz continuos with Lipschitz constant } 1.\vspace{0.1cm}\\
    \textbf{(ii) } \textit{ For any } v \notin \Lambda, \partial d_{\Lambda} (v) \textit{ has exactly one element with its length being } 1.\vspace{0.1cm}\\
    \textbf{(iii) } d_{\Lambda} (\cdot)^2 \textit{ is } C^1.
\end{array}
\end{equation*}
\label{AppendixMetric2Lemma}
\end{lemma}
Here, for a region $G \subset \R^n$ and a locally Lipschitz continuous function $\phi : G \rightarrow \R$ the \textit{Clarke's generalized gradient} $\partial \phi$ at $x \in G$ is defined as
\begin{equation*}
    \partial \phi (x) := \left\{ \xi \in \R^n : \langle \xi , y \rangle \leq \smash{\displaystyle\limsup_{z \rightarrow x, z \in G, t \downarrow 0}}\hspace{0.1cm} \frac{\phi (z + ty) - \phi(z)}{t} \right\},
\end{equation*}
which is nothing but the ordinary gradient $\nabla \phi (x)$ if $\phi$ is continuously differentiable at $x \in G$. Further, if $\phi : G \rightarrow \R$ is a convex function and the set $G \subseteq \R^l$ is convex it is a well-known fact that $\partial \phi (x)$ for $x \in G$ can be written as:
\begin{equation*}
    \partial \phi (x) = \{ \xi \in \R^l | \langle \xi , y \rangle \leq \phi(x+y) - \phi(x), \forall y \in \R^l, x + y \in G \}.
\end{equation*}
Thus, for convex functions defined over convex sets the generalized Clarke gradient coincides with the subgradient.

\end{document}